\documentclass[a4paper]{scrartcl}
\usepackage[colorlinks,bookmarksopen,bookmarksnumbered,citecolor=red,urlcolor=red]{hyperref}
\relax
\usepackage{todonotes}
\usepackage{stmaryrd}
\usepackage{moreverb}
\usepackage{graphicx}
\usepackage{amsthm}

\usepackage{amsmath,amssymb}
\usepackage[utf8]{inputenc}

\theoremstyle{definition}
\newtheorem{definition}{Definition}
\newtheorem{algorithm}{Algorithm}

\theoremstyle{plain}

\theoremstyle{remark}
\newtheorem{remark}{Remark}

\newcommand{\n}{\mathrm{n}} 

\newcommand{\fineGrid}{{\mathcal{T}_h}} 
\newcommand{\fineEl}{e} 
\newcommand{\fineIn}{f} 

\newcommand{\codimO}{\Gamma_a} 
\newcommand{\neuBounFP}{\Gamma_{p, n}} 
\newcommand{\dirBounFP}{\Gamma_{p, d}} 
\newcommand{\dirBounFS}{\Gamma_{s, d}} 
\newcommand{\innerInt}{\Gamma_i} 

\newcommand{\coarseGrid}{{\mathcal{T}_H}}
\newcommand{\coarseEl}{E}
\newcommand{\coarseIn}{F}
\newcommand{\codimOC}{\Gamma_A} 

\newcommand{\fineSpace}{{\mathcal{V}_h}} 
\newcommand{\fineSpaceVel}{{[\fineSpace]^d}} 
\newcommand{\coarseSpace}{{\mathcal{W}_N}} 
\newcommand{\hOne}{H^{1}} 
\newcommand{\lTwo}{L^{2}} 
\newcommand{\raviartThomas}{\text{RT}} 

\newcommand{\RR}{\mathbb{R}} 
\newcommand{\NN}{\mathbb{N}} 

\newcommand{\jump}[1]{\llbracket #1 \rrbracket}
\newcommand{\mean}[1]{\{#1\}}

\newcommand{\param}{\mathbf{\mu}} 
\newcommand{\params}{\mathcal{P}} 
\newcommand{\trainSet}{\mathcal{P}_\text{tr}} 
\newcommand{\nParam}{M}

\newcommand{\pBilForm}[3]{b_{h}(#1, #2; #3)}
\newcommand{\pLinForm}[3]{l_{h}(#1; #2, #3)}

\newcommand{\pdt}{\partial_t}

\newcommand{\dx}{\,\mathrm{d}x}
\newcommand{\ds}{\,\mathrm{d}S}

\newcommand{\dirValP}{{p_D}} 
\newcommand{\neuValP}{{v_N}} 
\newcommand{\dirValS}{{s_D}} 
\newcommand{\sInitial}{s_0}
\newcommand{\nwMob}{\gamma_o} 
\newcommand{\nwMobBC}{\lambda_o} 
\newcommand{\nwMobP}{\Lambda_o} 
\newcommand{\nwMobPP}[1]{\Lambda_{o}^{#1}} 
\newcommand{\wMob}{\gamma_w} 
\newcommand{\wMobBC}{\lambda_w} 
\newcommand{\wMobP}{\Lambda_w} 
\newcommand{\wMobPP}[1]{\Lambda_{w}^{#1}} 
\newcommand{\tMob}{\gamma} 
\newcommand{\tMobBC}{\lambda} 
\newcommand{\tMobPP}[1]{\Lambda^{#1}} 

\newcommand{\wnwMob}{\gamma_{\alpha}} 
\newcommand{\wnwMobBC}{\lambda_{\alpha}} 
\newcommand{\wnwMobPP}[1]{\Lambda_{\alpha}^{#1}} 

\newcommand{\nwVis}{\eta_o} 
\newcommand{\wVis}{\eta_w} 
\newcommand{\nwDensity}{\varrho_o} 
\newcommand{\wDensity}{\varrho_w} 
\newcommand{\perm}{\kappa} 
\newcommand{\por}{\phi} 
\newcommand{\f}{f_w} 
\newcommand{\numFlux}{\chi}
\newcommand{\gravity}{G}

\newcommand*{\p}[1][]{p^{#1}}
\newcommand*{\vel}[1][]{u^{#1}}
\newcommand*{\s}[1][]{s^{#1}}
\newcommand*{\tof}[1]{\tau^{#1}} 

\newcommand*{\pF}[1][]{{p_h^{#1}}} 
\newcommand*{\velF}[1][]{{u_h^{#1}}} 
\newcommand*{\sF}[1][]{{s_h^{#1}}} 
\newcommand*{\tofF}[1]{\tau^{#1}_h} 
\newcommand*{\sFS}[1][]{\tilde{s}_h^{#1}} 
\newcommand*{\pC}[1][]{{p_H^{#1}}} 
\newcommand*{\sC}[1][]{{s_H^{#1}}} 
\newcommand{\pD}[1][]{\underline{p}_H^{#1}}

\newcommand{\redBas}{\Phi}

\newcommand{\errorEst}{\Delta_N} 
\newcommand{\lerror}{\Delta_{\lTwo}} 
\newcommand{\herror}{\Delta_{\hOne}} 

\newcommand{\nTimeSteps}{N_T}
\newcommand{\tEnd}{t^{\nTimeSteps}}

\newcommand{\diam}{\text{diam}}
\newcommand{\ris}{\,\rule[-5pt]{0.4pt}{12pt}\,{}} 
\newcommand{\shockDe}{\mathcal{S}} 
\newcommand{\supp}[1]{\text{supp}(#1)} 
\newcommand{\vspan}[1]{\text{span}(#1)} 
\newcommand{\x}{x} 
\newcommand{\cf}{cf.\ } 
\newcommand{\transpose}{\tau}

\begin{document}
\title{The Localized Reduced Basis Multiscale method for two-phase flows in porous media}
\author{S.~Kaulmann\thanks{Institute for Computational and Applied Mathematics, Einsteinstr. 62, D-48149 Münster, Germany \url{{sven.kaulmann,mario.ohlberger}@uni-muenster.de}} \thanks{Institute for Applied Analysis and Numerical Simulation, Pfaffenwaldring 57, D-70569 Stuttgart, Germany, \url{{sven.kaulmann,haasdonk}@mathematik.uni-stuttgart.de}} ,
B.~Flemisch\thanks{Institute for Modelling Hydraulic and Environmental Systems, Pfaffenwaldring 61, D-70569 Stuttgart, Germany, \url{bernd@iws.uni-stuttgart.de}} , B.~Haasdonk$^\dagger$,\\ K.-A.~Lie\thanks{SINTEF ICT, PO Box 124. Blindern, N-0314 Oslo, Norway, \url{knut-andreas.lie@sintef.no}} , and M.~Ohlberger$^*$}
\date{}
\maketitle
\begin{abstract}
In this work, we propose a novel model order reduction approach for two-phase
flow in porous media by introducing a formulation in which the mobility, which
realizes the coupling between phase saturations and phase pressures, is regarded
as a parameter to the pressure equation. Using this formulation, we introduce
the Localized Reduced Basis Multiscale method to obtain a low-dimensional
surrogate of the high-dimensional pressure equation. By applying ideas from
model order reduction for parametrized partial differential equations, we are
able to split the computational effort for solving the pressure equation into a
costly \textit{offline} step that is performed only once and an inexpensive
\textit{online} step that is carried out in every time step of the two-phase
flow simulation, which is thereby largely accelerated. Usage of elements from
numerical multiscale methods allows us to displace the computational intensity
between the offline and online step to reach an ideal runtime at acceptable
error increase for the two-phase flow simulation.
\end{abstract}
%
\section{Introduction} \label{sec:introduction} Real-world applications of
two-phase flow in a porous medium are often characterized by large disparities
in physical scales in the sense that local variations in the permeability and
porosity of the porous medium may strongly affect flow patterns on scales that
are orders of magnitude larger. To accurately resolve the effects of local
variations, one easily ends up with equation systems (for pressures and phase
saturations) that may contain millions of unknowns. Although (non)linear systems
of this size obviously can be solved within a reasonable time-frame using
high-performance computing, there are also many engineering workflows in which
high-performance computing is infeasible or impossible. In this work, we will
therefore discuss methods for accelerating the computation.
To this end, we will focus on sequential simulation methods
in which the flow and transport are solved in separate steps, and in which
updating the flow pattern is the most computationally expensive part.

The classical approach to reducing the size of the discrete simulation model is
to upscale the problem. That is, use a numerical method for calculating
effective parameters and functions on a coarser scale that represent the local
flow effect of the unresolved scale in an averaged sense. Standard upscaling
methods solve representative fine-scale (flow or transport) problems to
calculate averaged quantities (such as permeabilities, transmissibilities,
etc.). There are both local and global upscaling methods. Local methods choose
sub-domains of a size much smaller than the global scale (e.g., the size of one
coarse grid block) and calculate effective parameters locally for each of these
sub-domains. Global methods solve representative fine-scale problems on the
global scale. Local methods can be further extended to (adaptive) local-global
methods. In this case, a downscaling step is added to the method to approximate
local fine-scale boundary conditions from the global coarse-scale solution.
Excellent overviews of upscaling method can be found in, e.g.,
\cite{Farmer2002,Durlofsky2005}.

A related, but more recent approach is to use so-called multiscale methods that
try to incorporate fine-scale information into the set of coarse-scale equations
in a way that is consistent with the local property of the differential operator
\cite{Hughes:1995:MPG,EfendievHou:09}. There are many different
multiscale methods: Hughes et al.\ \cite{HughesEtAl:1998:VMM} introduced the
Variational Multiscale method as a general framework for constructing
hierarchical approximation spaces for Dirichlet-type problems. The Multiscale
Finite-Element method \cite{HouWu:1997:MFE} uses solutions to local fine-scale
problems to incorporate fine-scale details into coarse-scale basis functions; a
recent development in this direction is the Generalized Multiscale
Finite-Element method \cite{Efendiev:2013wg}. The Mixed Multiscale
Finite-Element method \cite{ChenHou2003,Aarnes:2008jj} is based on the same
principles but in addition allows for mass-conserving reconstruction of the
fine-scale velocities. Related methods include the Numerical Subgrid Upscaling
method \cite{Arbogast:02} and the Multiscale Mortar method
\cite{ArbogastEtAl:2007:MMM}. In the Multiscale Finite-Volume method
\cite{JennyEtAl:2003}, the underlying idea is to construct coarse-scale
transmissibilities that account for fine-scale effects and lead to a multi-point
approximation for the Finite-Volume discretization on the coarse scale.
Recently, a more robust Two-Point Finite-Volume method was proposed
\cite{MoynerLie:14}. In \cite{Henning:2013tc}, an adaptive version of the
Heterogeneous Multiscale method \cite{EEngquist:03:HMM}, which aims at capturing
the macroscopic scale of the problem by estimating necessary information from a
microscopic model, was applied to immiscible two-phase flows in porous media.
Finally, Henning et al.\ \cite{Peterseim:2013vz} recently applied the Partition
of Unity method in the multiscale context as a means for reliable numerical
homogenization for elliptic equations with rough coefficients.

Reduced Basis (RB) methods represent an approach for model order reduction for
parametrized partial differential equations, which has gained large
popularity in the last decade. In numerous publications, these methods have proven to
be versatile tools for reducing the computational effort for stationary problems
\cite{Prudhomme:2002gm,Machiels:2000jx,Prudhomme:2002ks,Patera:un,Nguyen:2005gd}
as well as for linear \cite{Grepl:2005eu,Rovas:2006ke} and nonlinear
\cite{Grepl:2007kx} parabolic problems, and linear
\cite{Haasdonk:2008io,Haasdonk:2008ts} and nonlinear \cite{Haasdonk:2009uc}
hyperbolic problems. Applications to two-phase flow in porous media exist
\cite{Drohmann:2012tc}. Considered parametrizations cover initial, source and
boundary data, geometries and material- and control-factors.

The key idea of RB methods is to introduce two computational phases: an
offline and an online phase. During the computationally demanding offline phase,
a set of solutions to the parametrized equation at hand, so-called
\emph{snapshots}, are computed for different parameters using a Greedy-type
algorithm \cite{Veroy:2012cb}. From these snapshots, a reduced basis is computed
and the offline-phase is concluded by projecting the high-dimensional
discretization onto the reduced basis, giving rise to a reduced-dimensional
equation.

During the online phase, a solution for a given parameter is computed by solving
the reduced equation. Using the common assumption of
\emph{parameter-separability} or \emph{affine parameter dependence}
\cite{Patera:un}, a complete decoupling of the offline and online phases becomes
straightforward. Hence, the computational effort during the online phase does
not depend on the dimension of the underlying fine discretization and decisive
speedups over the latter are possible at marginal loss of accuracy. If the
assumption of parameter-separability cannot be fulfilled for the equation at
hand, techniques like the Empirical Interpolation method
\cite{Barrault:2004kz,Drohmann:2012ea} can be used.

Herein, we will study the Localized Reduced Basis Multiscale method (LRBMS),
which was originally introduced for general parametrized elliptic problems in
\cite{Kaulmann2011rr, Albrecht2012lq}, and which is not limited to the current
two-phase flow setting. The method allows parametrizations like
boundary, initial or source data, but those will be neglected in this work
in order to concentrate on the new parametrization introduced below. Application of the LRBMS method yields
a potentially costly offline phase, in which localized reduced-dimensional bases
for the pressure are computed. These bases then operate as low-dimensional
surrogates of the high-dimensional Discontinuous Galerkin discretization during
the actual two-phase flow simulation (here referred to as the online phase). The
online phase is carried out using a sequential splitting scheme to decouple the
high-dimensional saturation equation and reduced-dimensional pressure equation.

A key point in this contribution is the use of the so-called time-of-flight $\tau$ to
pa\-ra\-metrize our two-phase model. The time-of-flight is the time it takes for an inert
particle to travel along a pathline from the closest point on the inflow
boundary and to a given point inside the domain, or alternatively, the time it
takes for a particle to travel from a given interior point and to the closest point
on the outflow boundary.  In the absence of gravity, the three-dimensional Eulerian equations
that describe how fluid phases are transported, can be transformed into a family
of one-dimensional equations in Lagrangian coordinates (i.e., along streamlines). Here,
$\tau$ plays the role as the spatial coordinate for each one-dimensional transport equation.
The time-of-flight therefore carries valuable information about the
characteristics of the flow, and our idea is to use $\tof{s}(x)$, defined for a
given spatial point $x$ and a given saturation snapshot $s$, as a means to
efficiently compute a limited number of mobility profiles
$\tMobPP{1},\ldots,\tMobPP{\nParam}$, which are assumed to approximate the actual
mobility $\tMobBC(s)$ for a given saturation $s$ via $\tMobBC(s)\approx \sum
\mu_i \tMobPP{i}$ for a given parameter $\mu\in\RR^\nParam$. The coupling between
the saturation and pressure equation via the mobility is now replaced by a
coupling via the parameter $\mu$ and we apply model order reduction techniques
for parametrized problems to the pressure equation, now being parametrized by
$\mu$.

The approximation quality of the reduced system is typically controlled using
\emph{a-posteriori} error estimators.
In a recent work, Schindler et al.\ \cite{Schindler:2014tz} introduced a
localized \textit{a-posteriori} error estimate perfectly fit for the localized
method introduced in Section \ref{sec:lrbms}. Another suitable estimate was
recently introduced in \cite{Smetana:vx}.

The rest of the paper is structured as follows: After introducing the details of
the two-phase flow setting in the next section, we will describe its
high-dimensional discretization in Section~\ref{sec:highDimDisc}.
Section~\ref{sec:lrbms} is dedicated to the derivation of the Localized Reduced
Basis Multiscale method and Section~\ref{sec:couplingSatAndPres} gives details
on the coupling between the reduced pressure equation and the high-dimensional
saturation equation. Finally, in Section~\ref{sec:numericalExperiments}, we
demonstrate the applicability of the method in numerical experiments.

\section{Model Problem}
\label{sec:modelProb}
In this section, we introduce the mathematical model of two-phase flow that will
be used throughout the rest of this work. Our model problem is the flow of two
phases in a spatio-temporal domain $\Omega\times [0,T]\subset\RR^d\times\RR^+$
where $d\in\{2,3\}$ denotes the space dimension. We use a global pressure, total
velocity formulation for two incompressible, immiscible fluids that includes
gravity but no capillary effects. This yields the equation for the unknown
pressure $p$
\begin{align}
\label{eq:mainPressure}
-\nabla \cdot \Bigl(\tMobBC(s)\perm \nabla p - \perm \bigl[\wMobBC(s)\wDensity+\nwMobBC(s)\nwDensity\bigr]\gravity\Bigr) &= q_1 \quad \text{in } \Omega\times [0,T],
\end{align}
where $\wDensity$ and $\nwDensity$ denote the densities for the wetting and
non-wetting phases, respectively, and $G$ is the gravitational force vector
$\gravity = (0,0,-g)^\transpose$,  with $g$ being the gravitational
acceleration. Furthermore, $\perm$ denotes the total permeability and  the
functions $\wMobBC,\,\nwMobBC,\,\tMobBC:[0,1]\to\RR^+$ denote the wetting,
non-wetting and total mobility, given by
\begin{align}
\label{eq:brooksCorey}
{\wMobBC}(s) = \frac{k_\text{rw}(s)}{\wVis},\quad {\nwMobBC}(s)= \frac{k_\text{ro}(s)}{\nwVis}, \quad \tMobBC(s)=\wMobBC(s)+\nwMobBC(s),
\end{align}
where $\wVis$ and $\nwVis$ are the viscosities of the wetting and non-wetting
phase, respectively and the relative permeabilities $k_\text{rw}$ and
$k_\text{ro}$ are given functions of the  saturation $s$ of the wetting phase.

Using Darcy's law, the total velocity $u$ is given by
\begin{align}
\label{eq:mainVelocity}
u &= - \tMobBC(s) \perm \nabla p + \perm \bigl[\wMobBC(s)\wDensity + \nwMobBC(s)\nwDensity\bigr]\gravity \quad \text{in } \Omega\times [0,T],
\end{align}
and enters the transport equation for the saturation $s$,
\begin{align}
\label{eq:mainSaturation}
\phi\,\pdt s +\nabla\cdot\Bigl(f_w(s)\, \bigl[u+ \perm\nwMobBC(s)(\wDensity-\nwDensity)\gravity\bigr]\Bigr) &= q_2 \quad \text{in } \Omega\times [0,T].
\end{align}
Here, $\por$ denotes the porosity of the porous medium and the fractional flow
of water $\f$ is given by
\begin{align}
\label{eq:frationalFlow}
\f(s) = \frac{\wMobBC(s)}{\wMobBC(s)+\nwMobBC(s)}.
\end{align}
The above three equations for the pressure (\ref{eq:mainPressure}), velocity
(\ref{eq:mainVelocity}) and saturation (\ref{eq:mainSaturation}) are equipped
with the boundary conditions
\begin{align}
\begin{aligned}
\label{eq:mainBoundary}
s &= \dirValS &\text{in } \partial\Omega_{s,d}\times [0,T], \\
 - \tMobBC(s) \perm \nabla p\cdot\n &= \neuValP &\text{in } \partial\Omega_{p,n} \times [0,T],  \\
p &= \dirValP &\text{in } \partial\Omega_{p,d}\times [0,T],
\end{aligned}
\end{align}
on the inflow boundary for the saturation $\partial\Omega_{s,d}$ and on the
Dirichlet and Neumann boundaries for the pressure $\partial\Omega_{p,d}$,
$\partial\Omega_{p,n}$. Further, we impose initial conditions
\begin{align}
\label{eq:mainInitialData}
s(\cdot, 0) &= \sInitial \quad\text{in }\Omega.
\end{align}


\section{High-Dimensional Discretization}
\label{sec:highDimDisc}
This section is dedicated to the exposition of the so-called high-dimensional
discretization of problem (\ref{eq:mainPressure})--(\ref{eq:mainBoundary}). The
term high-dimensional is used here to indicate the antonym of the reduced
(low-dimensional) discretization, \cf Section~\ref{sec:introduction}.

Different schemes are frequently used to discretize Equations
(\ref{eq:mainPressure})-(\ref{eq:mainBoundary}). A tendency towards
finite-volume-type schemes can be observed in the literature because they
are motivated by local conservation properties. As it will prove advantageous in
the analysis
of our reduced method, we will use a Discontinuous Galerkin (DG) discretization with
arbitrary local polynomial degree. In particular, we have found the
Symmetric Weighted Interior-Penalty (SWIP) DG method \cite{Ern:2010cn},
including the total velocity reconstruction presented therein, to be very robust
and it is therefore our method of choice in the following.

\subsection{Discretization}
As a first step towards a discrete version of Equations
(\ref{eq:mainPressure})--(\ref{eq:mainBoundary}), we introduce an admissible
tessellation $\fineGrid$ of the computational domain $\Omega$. We assume
$\fineGrid$ to be a set of non-overlapping elements $\fineEl$ with width
$h_\fineEl = \diam(\fineEl)$ and define the grid size $h =
\max_{\fineEl\in\fineGrid} h_\fineEl$. The intersections of co-dimension one of
an element $\fineEl$ are denoted by $\fineIn$ and their width by $h_\fineIn =
\diam(\fineIn)$. With each intersection $\fineIn$ we associate a unique normal
vector $\n_\fineIn$, for which the subscript $\fineIn$ will be dropped whenever
no ambiguity arises. We let $\dirBounFP$ and $\neuBounFP$ denote all
boundary intersections of $\fineGrid$ where Dirichlet and
Neumann conditions for the pressure shall be implied, respectively, and
$\dirBounFS$ denote all boundary intersections on which we impose a fixed
saturation. Finally, by $\innerInt$ we denote the inner
intersections and by $\codimO =
\neuBounFP\cup\dirBounFP\cup\innerInt$ the set of all
intersections.

Additionally we introduce a temporal
discretization $t^0,\ldots,\tEnd\in[0,T]$ and space-time-discrete approximations
$\sF[n]\approx s(\cdot, t^n)$, $\pF[n]\approx p(\cdot, t^n)$ of the saturation
and the pressure, both stemming from the space of piecewise polynomials of
maximum local order $k$
\begin{align}
\fineSpace = \left
				\{v\in \lTwo(\Omega) \left|
					v_{|\fineEl}\in\mathbb{P}^k(\fineEl)\quad \forall\fineEl\in\fineGrid
				\right.
			\right\}.
\end{align}
Functions in $\fineSpace$ are two-valued on intersections $\fineIn =
\fineEl_1\cap\fineEl_2$ (with $\n_\fineIn$ pointing from $\fineEl_1$ to
$\fineEl_2$). We define jump $\jump{\cdot}_f$ and (weighted) mean $\mean{\cdot}_f$ values
for $w\in\fineSpace$, i.e.,
\begin{align*}
\jump{w}_f &= w_{|\fineEl_1}-w_{|\fineEl_2},\\
\mean{w}_f &= \tau_{\fineEl_1,\fineIn} w_{|\fineEl_1} + \tau_{\fineEl_2,\fineIn} w_{|\fineEl_2},
\quad \tau_{\fineEl_{\ell},\fineIn} = \frac{a_{\ell,\fineIn}}{a_{1,\fineIn}+a_{2,\fineIn}},
\quad a_{\ell, \fineIn}=\|(\perm)_{|\fineEl_{\ell}}\|_{L^\infty(\fineIn)},\; \ell=1,2.
\end{align*}
For the sake of readability, the subscript $\fineIn$ in $\jump{\cdot}_\fineIn$
and $\mean{\cdot}_\fineIn$ will be dropped in the following whenever no
ambiguity arises.

\subsection{Pressure Equation} Following the ideas in \cite{Ern:2010cn}, we
introduce discrete formulations of Equations
(\ref{eq:mainPressure})--(\ref{eq:mainBoundary}). For a given general mobility
function $\tMob:\Omega\to\RR^+$ the bilinear form
$\pBilForm{\cdot}{\cdot}{\tMob}:\fineSpace\times\fineSpace\to\RR$ is defined as
\begin{align}
\label{eq:swipBilP}
\begin{split}
\pBilForm{v}{w}{\tMob} &= \sum_{\fineEl\in\fineGrid}\int_\fineEl \tMob \perm \nabla v \nabla w \dx
 + \sum_{\fineIn\in\innerInt\cup\dirBounFP} \frac{\sigma_\fineIn}{h_\fineIn} \int_\fineIn \jump{v}\jump{w} \ds \\
 &\quad + \sum_{\fineIn\in\innerInt\cup\dirBounFP}\int_\fineIn
 \left( \mean{\tMob\perm\nabla v\cdot \n_\fineIn}\jump{w} +  \mean{\tMob\perm\nabla w\cdot \n_\fineIn}\jump{v} \right) \ds.
\end{split}
\end{align}
Here, the penalty parameter $\sigma_\fineIn$ is given as
\begin{align}
\label{eq:penaltyParam}
\sigma_\fineIn = c_\fineIn\frac{2a_{1,\fineIn}a_{2,\fineIn}}{a_{1,\fineIn}+a_{2,\fineIn}},
\end{align}
where $c_\fineIn>0$ is a constant that has to be chosen larger than a
minimal threshold depending on the regularity of the mesh.

\begin{remark}
Note that in \cite{Ern:2010cn}, $a_{\ell,
\fineIn}=\|(\tMob\perm)_{|\fineEl_{\ell}}\|_{L^\infty(\fineIn)}$ was used. As the
penalty parameter $\sigma_\fineIn$ would not allow an affine decomposition, we
refrain from using the mobility in the penalties and instead replace the
original constant $c_\fineIn$ by ${c_\fineIn}/{\max\{\wVis, \nwVis\}}$. Usage of
the original weighted averages and penalties would be possible using the
Empirical Interpolation technique \cite{Barrault:2004kz,Drohmann:2012ea}.
\end{remark}

The right-hand side for the discrete formulation of the pressure equation is for
$s\in\fineSpace$, $\nwMob,\wMob:\Omega\to\RR^+,\,\tMob=\nwMob+\wMob$ given via
the linear form $\pLinForm{\cdot}{\nwMob}{\wMob}:\fineSpace\to\RR$,
\begin{align}
\begin{split}
\pLinForm{w}{\nwMob}{\wMob} &=
	\sum_{\fineEl\in\fineGrid}\int_{\fineEl} q_1 w + (\nwMob\nwDensity+\wMob\wDensity)\perm\,G\cdot\nabla w\dx \\
	&\quad- \sum_{\fineIn\in\innerInt\cup\dirBounFP}\int_\fineIn \mean{(\nwMob\nwDensity+\wMob\wDensity)\perm\,G\cdot\n_\fineIn}\jump{w} \ds \\
	&\quad+ \sum_{\fineIn\in\dirBounFP}\int_{\fineIn} \left( \frac{\sigma_\fineIn}{h_\fineIn} w-\tMob\perm\nabla w\cdot \n_\fineIn \right) \dirValP \ds
	+ \sum_{\fineIn\in\neuBounFP}\int_\fineIn \neuValP w \ds.
\end{split}
\end{align}
At this point we are ready to introduce the so-called high-dimensional
discretization of (\ref{eq:mainPressure}).

\begin{definition}
\label{def:pressureDG}
For given functions $\wMob,\nwMob:\Omega\to\RR^+$ and
$\tMob=\wMob+\nwMob$, we call $\pF\in\fineSpace$ the high-dimensional pressure
solution if it satisfies
\begin{align*}
\pBilForm{\pF}{w}{\tMob} = \pLinForm{w}{\nwMob}{\wMob}\quad\forall w\in\fineSpace.
\end{align*}
\end{definition}

\subsection{Total Velocity Reconstruction}
From a given pressure $\p$, which may be either $\pF$ or a reduced pressure later on, we compute the total velocity $\velF$ using the
conservative reconstruction from \cite{Ern:2010cn}, which yields continuity of
the normal components. For $k=1$, i.e., piecewise linear pressure, the velocity
is computed in the lowest order Raviart-Thomas space
\begin{align}
\raviartThomas = \left\{ u\in H(\text{div})	\left| u\ris_T\in [\mathbb{P}_0(T)]^d+\mathbf{x}\mathbb{P}_0(T)\; \forall T\in\fineGrid \right.\right\},
\end{align}
where $H(\text{div}) = \{v\in [\lTwo(\Omega)]^d \left| \nabla\cdot v\in
\lTwo(\Omega)\right. \}$ and $\mathbf{x}=(x_1,\ldots,x_d)^\tau$. Additionally, we
introduce
\begin{align*}
\jump{v}_\fineIn^* = \begin{cases}
\jump{v}_\fineIn, & \fineIn\in\innerInt, \\
v\ris_\Omega-v_D, & \fineIn\in\Gamma_{v,d},\\
0, & \fineIn\in\neuBounFP,
\end{cases}
\end{align*}
for $v\in\fineSpace$ with Dirichlet boundary data $v\ris_{\Gamma_{v,d}} = v_D$.
\begin{definition}
\label{def:velocityDG}
For a given pressure $p\in\fineSpace$ and saturation
$s\in\fineSpace$, let the total velocity $\velF\in\text{RT}$ be the unique
solution to
\begin{align}
\label{eq:velocityDG}
\int_{\fineIn} (\velF\cdot\n)\ds
	= \int_\fineIn -\n\cdot\mean{\tMobBC(\s)\perm\nabla\p-\perm \bigl[\wMobBC(s)\wDensity+\nwMobBC(s)\nwDensity\bigr]\gravity}
	+ \frac{\sigma_\fineIn}{h_\fineIn}\jump{\p}_\fineIn^* \ds\quad\forall \fineIn\in\codimO.
\end{align}
\end{definition}

\subsection{Saturation Equation}
The saturation equation discretized by a DG scheme in space
using an upwind flux in the space of piecewise polynomial functions
$\fineSpace$, and for the temporal discretization, we use an explicit Euler
scheme.
\begin{definition}
\label{def:saturationDG}
For a given saturation $\sF[n]\in\fineSpace$ and  velocity
$\vel\in\fineSpaceVel$, the saturation $\sF[n+1]\in\fineSpace$ is given as the
solution to %
\begin{align}
\label{eq:saturationDG}
\begin{split}
\sum_{\fineEl\in\fineGrid}\int_\fineEl \frac{\por}{\Delta t}\sF[n+1] v_h \dx &=
	\sum_{\fineEl\in\fineGrid}\int_\fineEl \frac{\por}{\Delta t}\sF[n] v_h \dx
	+ \sum_{\fineEl\in\fineGrid} q_2 v_h\dx \\
	&\quad +\sum_{\fineEl\in\fineGrid} \int_\fineEl \bigl[\vel+\nwMobBC(\sF[n])(\wDensity-\nwDensity)\perm\cdot\gravity\bigr] \f(\sF[n])  \cdot \nabla v_h\dx \\
	&\quad - \sum_{\fineIn\in\codimO} \int_\fineIn \n_\fineIn \cdot \mean{\vel+\nwMobBC(\sF[n])(\wDensity-\nwDensity)\perm\cdot\gravity}\, \f(\numFlux(\sF[n]))\jump{v_h}\ds \\
	&\quad -\sum_{\fineIn\in\innerInt\cup\dirBounFS} \frac{\sigma_\fineIn}{h_\fineIn}\int_\fineIn \jump{\sF[n]}^*_f\jump{v_h}_f\ds,
\end{split}
\end{align}
for all $v_h\in\fineSpace$. Here the upwind function $\numFlux$ is given as
\begin{align}
\label{eq:upwind}
\numFlux(\sF[n]) =
	\begin{cases}
		{\sF[n]}^\uparrow, & \fineIn\in\innerInt, \\
		{\dirValS}, & \fineIn\in\dirBounFS, \\
	\end{cases}
\end{align}
where ${\sF[n]}^\uparrow$ denotes the upwind value of $\sF[n]$: If
$\fineIn\in\innerInt$, there exist $\fineEl_1, \fineEl_2\in\fineGrid$ such that
$\fineIn=\partial\fineEl_1\cap\partial\fineEl_2$. Then
${\sF[n]}^\uparrow={\sF[n]}\ris_{\fineEl_1}$ if $\velF \cdot \n_\fineIn\ge 0$
and ${\sF[n]}^\uparrow={\sF[n]}\ris_{\fineEl_2}$ otherwise. Remember that we
assume $\n_\fineIn$ to point from $\fineEl_1$ to $\fineEl_2$. The penalty
$\sigma_\fineIn$ in Equation (\ref{eq:saturationDG}) is the same as for the
pressure equation, see (\ref{eq:penaltyParam}).
\end{definition}

\subsection{Slope Limiter}
Discontinuous Galerkin methods require some kind of
stabilization to avoid over- and undershoots if they are to be used in a
two-phase flow context; see \cite{Riviere:2008vf,Dedner:2011wu}, for example.
Both artificial diffusion and slope limiters have been used as a means to this
end. We will employ a slope limiter introduced by Dedner et al.\
\cite{Dedner:2011wu} that is applicable to different element types and in two
and three spatial dimensions. We will shortly outline the most important
features of the mass-conservative limiter that we use in our experiments for the
case of a piecewise linear saturation, where we use the ``$\mathcal{DG}$
scheme'' (in contrast to ``$\mathcal{DG+R}$'', see
\cite[Section~6.1]{Dedner:2011wu}). For more details on the method and a
discussion of the general case of arbitrary polynomial order on $\sF[n]$ we
refer to \cite{Dedner:2011wu}.

The first step is to compute a so-called \emph{shock-detector}
$\shockDe:\fineGrid\to\RR^+$, which for a given saturation $\sF\in\fineSpace$
and velocity $\velF\in\fineSpaceVel$ reads:
\begin{align}
\shockDe(\fineEl) = \sum_{\mathcal{I}_\fineEl^\uparrow}\left(0.08\cdot d\cdot\sqrt{h_\fineEl}|\fineEl|\right)^{-1}\int_{\fineIn} \jump{\sF}\ds,
\end{align}
where $\mathcal{I}_\fineEl^\uparrow = \left\{\fineIn\in\codimO\left|
\fineIn\subset\partial\fineEl, \velF\cdot\n_\fineIn<0 \right.\right\}$ denotes
the upstream interfaces of $\fineEl$. Using $\shockDe$ makes it possible to
apply the limiter (i.e., reduce gradients of the saturation) only on cells in
which (strong) discontinuities are present and therefore unwanted numerical
oscillations may occur. In other regions, the saturation is left unlimited.
Cells that will be flagged in the above sense are then all cells
$\fineEl\in\fineGrid$ with $\shockDe(\fineEl)>1$ or $\sF(x)\notin [0,1]$ for
some $x\in\fineEl$. Let $\fineEl^*\in\fineGrid$ be such a cell and $\{\fineEl_i
\left|\, i=1,\ldots,N_n\right.\}$ its direct neighbors. We then compute
\begin{align*}
g_i = \nabla\sF\cdot(b_{\fineEl_i}-b_{\fineEl^*})\quad\text{and}\quad
d_i = \frac{1}{|\fineEl_i|}\int_{\fineEl_i} \sF\dx - \frac{1}{|\fineEl^*|}\int_{\fineEl^*} \sF\dx,
\end{align*}
where $b_\fineEl$ denotes the barycenter of the cell $\fineEl$. From $g_i$ and
$d_i$ we compute the gradient scales $m_i$ by
\begin{align*}
m_i =
	\begin{cases}
	0, & \text{if } g_i d_i<0, |g_i|>10^{-8} \text{ and } |d_i|>10^{-8}, \\
	d_i/g_i, & \text{if } g_i d_i>0,|g_i|>|d_i|, |g_i|>10^{-8}\text{ and } |d_i|>10^{-8}, \\
	1, &\text{otherwise}.
	\end{cases}
\end{align*}
Finally, the stabilized saturation $\sFS$ is computed as
\begin{align}
\label{eq:limiterMain}
\int_\fineEl \sFS\varphi\dx = \int_\fineEl \sF\varphi\dx+\int_\fineEl\min_i(m_i) \nabla\sF\cdot(x-b_\fineEl)\varphi\dx\quad\forall \varphi\in\mathbb{P}_1(\fineEl)
\end{align}
in all flagged cells $\fineEl\in\fineGrid$. In all other cells, we set $\sFS=\sF$.

\begin{algorithm}[High-Dimensional Two-Phase Flow Scheme]
\label{alg:highDimTwoPhaseFlow}
The high-dimensional simulation scheme for the two-phase flow is as follows:
\begin{enumerate}
\item Project the initial data $\sInitial$ onto the high-dimensional discrete
space: Let $\sF[0]\in\fineSpace$ be given by
	\begin{align*}
		\int_{\Omega} \sInitial \psi\dx = \int_\Omega \sF[0]\psi\dx\quad\forall \psi\in\fineSpace.
	\end{align*}
\item For all $n\in\{0,\ldots,\nTimeSteps-1\}$,
	\begin{enumerate}
	\item use the pressure scheme (Definition \ref{def:pressureDG}) for
	$\nwMob=\nwMobBC(\sF[n])$, $\wMob=\wMobBC(\sF[n])$ to compute $\pF[n+1]$;
	\item use the velocity scheme (Definition \ref{def:velocityDG}) for
	$p=\pF[n+1]$ to compute the total velocity $\velF[n+1]$;
	\item use the saturation scheme (Definition \ref{def:saturationDG}) for
	$u=\velF[n+1]$ to compute the saturation $\overline{s}_h^{n+1}$;
	\item use the limiter (Equation (\ref{eq:limiterMain})) for
	$\sF=\overline{s}^{n+1}_h$, $\velF=\velF[n+1]$ to compute a stabilized version
	$\sFS[n+1]$ of the saturation and set $\sF[n+1]=\sFS[n+1]$.
	\end{enumerate}
\end{enumerate}
\end{algorithm}

\subsection{Time-of-Flight Equation}
\label{sec:timeOfFlight} The so-called
time-of-flight $\tof{}(\x)$ is defined as the time it takes for a passive
particle to reach a given point $\x\in\Omega$, starting from the closest point
on the inflow boundary. Here, $\tof{}(\x)$ can be defined by integrating $\int
\por/|\vel{}|\,dx$ along streamlines, or by solving $\vel{}\cdot\tof{}=\por$.
The time-of-flight inherits important information about the flow pattern and
will be used here to approximate the spatio-temporal behavior of the saturation
without computing a whole temporal evolution the transport equation.

Consistent with our discretization of the saturation and pressure equation we
use the DG-discretization introduced in \cite{Natvig:2007tr}: Find
$\tofF{s}\in\fineSpace$ such that
\begin{align}
\label{eq:timeOfFlight}
-\int_\fineEl (\tofF{s} \velF)\cdot\nabla\psi_h\dx + \int_{\partial\fineEl} ({\tofF{s}}^\uparrow\velF\cdot\n)\psi_h\ds
	=  \int_\fineEl \por \psi_h\dx,\quad\forall\psi_h\in\fineSpace
\end{align}
for all $\fineEl\in\fineSpace$. Here, ${\tofF{s}}^\uparrow$ again denotes the
upwind value of the two-valued function $\tofF{s}$, see (\ref{eq:upwind}), and
$\por$ denotes the porosity of the porous medium.

In the absence of gravitational forces, the discretized time-of-flight equation can be
permuted to a lower block-triangular form---if the computational mesh is reordered
according to the direction of the flow---and hence solved very efficiently in a
per-element fashion by a simple backsubstitution method; see
\cite{Natvig:2007tr} for details.


\section{The Localized Reduced Basis Multiscale Method}
\label{sec:lrbms}

Reduced Basis (RB) methods as introduced in Section~\ref{sec:introduction} offer
good complexity reduction for a lot of different kinds of equations.
Nevertheless, for real-world multiscale applications, the number of unknowns in
the high-dimensional discretization can be prohibitively large such that the
Greedy-type algorithm \cite{Veroy:2012cb} becomes unfeasible
because the number of snapshots computed during the basis generation cannot be
controlled.

A range of different extensions to the original RB method that can be seen as a
remedy for this problem were developed in recent years. The Reduced Basis
Element method \cite{Knezevic:2012cf,Maday:2004iw,Chen:2010fb} introduces a
number of representative domains and computes reduced bases for each one of
them. The domain under consideration is then built up from the representative
domains (possibly using deformation mappings) and computations are carried out
in a space built up from the precomputed spaces. The Reduced Basis Hybrid method
\cite{Iapichino:2012go} and the Reduced Basis Domain Decomposition Finite
Element method \cite{Iapichino:2012vk} make use of this concept. In
\cite{Maier:vc}, Maier and co-workers introduce a RB approach for
heterogeneous domain decomposition that could be used in the current setting to
reduce the number of global snapshots, too.
Another approach for model order reduction for our application is presented in
\cite{Chaturantabut:2011cl}, in which Discrete Empirical Interpolation and
Proper Orthogonal Decomposition are used in conjunction to build reduced-order
models for simulation of viscous fingering in porous media.

In our recent articles \cite{Albrecht2012lq, Kaulmann2011rr}, we introduced the
so-called Localized Reduced Basis Multiscale method (LRBMS). It connects ideas
from numerical multiscale methods with the RB approach and aims at reducing the
offline time of the RB method while ensuring at least identical approximation
quality during the online time, possibly at a slightly increased cost in terms
of runtime. (Similarly, in \cite{Krogstad:11:POD}, the proper orothogonal
decomposition method is combined with a multiscale mixed finite element
framework.) The main idea of this approach is to introduce two grids: A fine
triangulation for the high-dimensional discretization and a coarse one that is
used to build up localized RB spaces. The rest of this section is dedicated to
the exposition of the details of this method.

\subsection{Parameterization}
\label{sec:lrbmsParameterization}
At this point we provide details on the parametrization indicated in
Section~\ref{sec:introduction}. As stated earlier, the coupling between the
saturation and pressure equation will be realized via a parametrization of the
pressure equation: We introduce wetting and non-wetting phase
mobilities $\wMobP,\nwMobP:\Omega\to\RR^+$ that depend on the saturation $s$ via
\begin{align}
\label{eq:parameterizedMob}
\wMobP(s)(\x) = \sum_{q=1}^{\nParam} \theta_q(s) \wMobPP{q}(\x),\quad
\nwMobP(s)(\x) = \sum_{q=1}^{\nParam} \theta_q(s) \nwMobPP{q}(\x),
\quad \x\in\Omega,
\end{align}
where $\wMobPP{q},\,\nwMobPP{q}\in\fineSpace$, $q=1,\ldots,\nParam$ denote
saturation-independent mobility profiles that were precomputed during the
offline phase and the coefficients $\theta_q$ are computed during the online
phase via a least-squares fitting: Let $\Theta(s)\in\RR^{\nParam}$ be given by
\begin{align}
\label{eq:leastSquaresFit}
\Theta(s) = \text{arg}\min_{\vartheta\in\RR^M }\|\tMobBC(s)-\sum_{q=1}^{\nParam}\vartheta_j\tMobPP{q} \|_{l^2},
\end{align}
then $\theta_q(s) = (\Theta(s))_q$. Here we used
$\tMobPP{q}=\wMobPP{q}+\nwMobPP{q}$. Details on the computation of the mobility
profiles $\nwMobPP{q}$, $\wMobPP{q}$ are given below.

\subsection{Discretization}
As mentioned before, we introduce a second, coarser
triangulation of the domain $\Omega$: Let $\coarseGrid$ be an admissible
triangulation with grid size $H = \max_{\coarseEl\in\coarseGrid}H_\coarseEl =
\max_{\coarseEl\in\coarseGrid}\text{diam}(\coarseEl)$, for  $H\gg h$. Cells in
$\coarseGrid$ will be denoted by $\coarseEl$, intersections of two elements by
$\coarseIn$, and the set of all coarse intersections by $\codimOC$.
Furthermore, we assume the two triangulations to be matching in the sense that
for each $\coarseIn\in\codimOC$ there exist $m\in\NN$ and
$\fineIn_1,\ldots,\fineIn_m\in\codimO$ such that $\coarseIn = \bigcup_{i=1}^m
\fineIn_i$.

Based on the coarse triangulation $\coarseGrid$, we introduce a localized
reduced-dimensional function space $\coarseSpace$. For every coarse grid cell
$\coarseEl\in\coarseGrid$ we assume a set of linearly independent functions
$\varphi_1^\coarseEl, \ldots, \varphi_{N^\coarseEl}^\coarseEl\in\fineSpace$ with
$\supp{\varphi_i^\coarseEl}\subset\coarseEl$ to be given so that they form the
local reduced basis: $\redBas^\coarseEl = \{ \varphi_1^\coarseEl,
\ldots, \varphi_{N^\coarseEl}^\coarseEl \}$. From those local bases we build a
global reduced basis $\redBas_N$ of size
$N=\sum_{\coarseEl\in\coarseGrid}N^\coarseEl$ on the whole domain $\Omega$ by
setting
\begin{align}
\redBas_N = \bigcup_{\coarseEl\in\coarseGrid} \redBas^\coarseEl.
\end{align}
We call the space $\coarseSpace=\vspan{\redBas_N}$, spanned by the global
reduced basis, \textit{reduced broken space} and point out that
$\coarseSpace\subset\fineSpace$ because of our DG ansatz.

Using the reduced basis $\redBas_N$, we compute the coarse-scale pressure
approximation $\pC\in\coarseSpace$.
\begin{definition}
\label{def:coarseScalePressure}
For a given reduced broken space $\coarseSpace$
and mobilities $\nwMob$, $\wMob,\tMob:\Omega\to\RR^+$ where
$\tMob=\nwMob+\wMob$, the reduced solution $\pC\in\coarseSpace$ to Equation
(\ref{eq:mainPressure}) is given by
\begin{align}
\label{eq:coarseScalePressure}
\pBilForm{\pC}{w}{\tMob} &= \pLinForm{w}{\nwMob}{\wMob},\quad\forall w\in\coarseSpace.
\end{align}
\end{definition}

\subsection{Basis Construction}
\label{sec:lrbmsBasisConstruction}
One of the key ingredients of the method proposed in this paper is a basis
construction algorithm that is a modification of the Greedy-type algorithm used
in conventional RB methods as mentioned in Section~\ref{sec:introduction}.
\begin{algorithm}[LRBMS Basis Construction]
\label{alg:lrbmsBasisConstruction}
Based on the coarse triangulation
$\coarseGrid$ and an error measure $\errorEst$, we compute the global reduced basis $\redBas_N$ using the
following algorithm:
\begin{enumerate}
\item Initialization
\begin{enumerate}
\item \label{itm:lrbmsBasisConsInitStart} Given the initial saturation
$\sInitial$, compute the initial pressure $\pF[0]\in\fineSpace$ and velocity
$\velF[0]\in\raviartThomas$.
\item From $\velF[0]$ compute the time-of-flight $\tofF{\sInitial}$.
\item \label{itm:lrbmsBasisConsInitEnd} From the time-of-flight and
$\nParam\in\NN$, $\nParam\ge 1$, compute mobility profiles
$\nwMobPP{1},\ldots,\nwMobPP{\nParam}\in\fineSpace$ and
$\wMobPP{1},\ldots,\wMobPP{\nParam}\in\fineSpace$ via
	\begin{align*}
		\wnwMobPP{1}(\x) &= \wnwMobBC(0) \\
		\wnwMobPP{\nParam}(\x) &= \wnwMobBC(1) \\
		\wnwMobPP{q}(\x) &= \begin{cases}
			\wnwMobBC(0), & \text{if } \tofF{\sInitial}(\x)>(q-1)\cdot\frac{T}{\nParam-2}, \\
			\wnwMobBC(1), & \text{otherwise},
		\end{cases}\quad \forall q\in{2,\ldots,\nParam-1}
	\end{align*}
	where $\wnwMobBC$ denotes the linear mobility (\ref{eq:brooksCorey}) of phase $\alpha$. Set $\tMobPP{q}=\wMobPP{q}+\nwMobPP{q}$ for $q=1,\ldots,\nParam$.
\item Choose a desired maximum basis size $N_{\max}\in\NN$, an approximation
tolerance $\epsilon_\text{tol}$, a parameter set $\params = \left\{
\param\in[\RR^+]^\nParam | \sum_i \param_i \le 1  \right\}$ and a discrete
subset $\trainSet\subset\params$, the so-called training set. Furthermore let
$\Phi^\coarseEl=\emptyset$ for all $\coarseEl\in\coarseGrid$,
$\Xi=\bigcup_{\coarseEl\in\coarseGrid}\Phi^\coarseEl$, $W=\vspan{\Xi}$ and
$\Psi=\emptyset$.
\end{enumerate}
\item Basis Extension
\begin{enumerate}
\item \label{itm:lrbmsGreedyStart}For each $\param\in\trainSet$ compute
$\pC(\tMob)\in W$ as the solution to the coarse-scale pressure equation
(\ref{eq:coarseScalePressure}), that is
\begin{align*}
\pBilForm{\pC(\tMob)}{w}{\tMob} = \pLinForm{w}{\nwMob}{\wMob}\quad\forall w\in W
\end{align*}
with $\wnwMob=\sum_{q=1}^{\nParam} \param_q \wnwMobPP{q}$ and
$\tMob=\wMob+\nwMob$.
\item Evaluate the error estimator for each parameter in the training set:
$\epsilon_\param:=\errorEst(\param, W)$ and find the parameter worst
approximated in the current basis: $\param_{\max} =
\text{arg}\max_{\param\in\trainSet}\epsilon_\param$.
\item\label{itm:lrbmsGreedySnapshot} If
$\epsilon_{\param_{\max}}>\epsilon_\text{tol}$ compute $\pF({\tMob^{\max}})$ as
solution to the high-dimensional equation (see Definition \ref{def:pressureDG})
for the mobilities $\wMob^{\max}=\sum (\param_{\max})_q\, \wMobPP{q}$ and
$\nwMob^{\max}=\sum (\param_{\max})_q\, \nwMobPP{q}$. Otherwise go to step
\ref{itm:lrbmsGreedyFinish}.
\item Extend the set of snapshots: $\Psi = \Psi\cup \left\{\pF({\tMob^{\max}})\right\}$.
\item For each $\coarseEl\in\coarseGrid$:
\begin{enumerate}
\item Use the Gram-Schmidt algorithm to orthonormalize the restriction of the pressure snapshot
$\pF({\tMob^{\max}})\ris_\coarseEl$ with respect to the current local basis
$\Phi^\coarseEl$:
\begin{align*}
	\widetilde p = \text{ortho}(\pF({\tMob^{\max}})\ris_\coarseEl, \Phi^\coarseEl).
\end{align*}
\item Extend the local reduced basis: $\Phi^\coarseEl =
\Phi^\coarseEl\cup\left\{\widetilde p\right\}$.
\end{enumerate}
\item Set $\Xi=\bigcup_{\coarseEl\in\coarseGrid}\Phi^\coarseEl$ and $W=\vspan{\Xi}$.
\item If $|\Xi|< N_{\max}$, go back to step \ref{itm:lrbmsGreedyStart}.
\end{enumerate}
\item Data Compression (optional)
\begin{enumerate}
\item \label{itm:lrbmsGreedyFinish} Set $\widetilde{\redBas}^\coarseEl =
\Psi\ris_\coarseEl$, that is: Define local bases $\widetilde{\redBas}^\coarseEl$
per coarse element $\coarseEl\in\coarseGrid$ as restrictions of the snapshots
$\Psi$ to the coarse elements.
\item\label{itm:lrbmsGreedyPCA} Apply the principal component analysis with
tolerance $\epsilon_\text{PCA}$ on each element:
\begin{align*}
\redBas^\coarseEl = \text{PCA}(\widetilde{\redBas}^\coarseEl, \epsilon_\text{PCA}).
\end{align*}
The details of this step can be found in the literature, see
\cite{Jolliffe:2002wx}, for example.
\end{enumerate}
\item Finalization
\begin{enumerate}
\item Define the global reduced basis $\redBas_N$ for
$N=\sum_{\coarseEl\in\coarseGrid} |\redBas^\coarseEl|$ as
\begin{align*}
\redBas_N = \bigcup_{\coarseEl\in\coarseGrid} \redBas^\coarseEl.
\end{align*}
\item Define the global reduced basis space $\coarseSpace$ as
\begin{align*}
\coarseSpace = \vspan{\redBas_N}.
\end{align*}
\end{enumerate}
\end{enumerate}
\end{algorithm}

In summary, the basis construction performs the following: In an
``initialization'' step, we compute pressure and velocity from the initial
saturation data. From the velocity, we compute approximate mobility profiles
using the time-of-flight. This will incorporate important features of the
problem, like low-permeability-lenses, for example, into the mobility profiles
and therefore also into the reduced basis for the pressure. After fixing some
input data like the desired basis size, we proceed to the basis extension step.

In the basis extension step, we add localized orthonormalizations of
high-dimensional snapshots to an initially empty local bases until either the
maximum total basis size is reached or a prescribed error tolerance, measured by
an \emph{a-posteriori} error estimator, is fulfilled in the current overall
basis, which is the joint of all local bases. Additionally, we save the original
snapshots.

In an optional ``compression'' step we define local per-coarse-element bases by
restricting the global untouched snapshots from the last step to each coarse
element. Next, we apply a data compression algorithm, the principal component
analysis, to each local basis to reduce the basis size in regions where
redundant information may be present, like in regions far from sinks and
sources, for example. At this point, it is important that we did not use the
orthonormalized local bases from the extension step for the local compressions.
Redundancies were canceled out in those bases and therefore the data compression
would not give meaningful results.

We conclude the algorithm with the ``finalization'' step by defining one global
reduced basis as the joint of all local bases and the reduced broken space as
its linear span.

The main ideas of this novel method are the restriction of the basis to elements
of a coarser grid and subsequent per-element data compression. By the
restriction to a coarse grid we reduce the number of snapshots needed to fulfill
a desired error tolerance on a prescribed training set of parameters during the
offline phase. This is easy to see for the limit case in which
$\coarseGrid=\fineGrid$, as in this case the reduced space coincides with the
high-dimensional discrete function space after a finite number of basis
extensions. For $|\coarseGrid|\le|\fineGrid|$, this effect was demonstrated in
\cite{Albrecht2012lq} and can be seen again in the numerical experiments in
Section~\ref{sec:numericalExperiments}.

The downside of this approach is that the size of the global reduced basis
increases with the number of coarse elements. This is where the idea of
per-element data compression comes into play and allows us to keep the total
basis size $N$, which is the main factor in online computation complexity, in
an agreeable range by reducing the local basis size in regions where little or no
variation is inherent in the local bases, as may be the case in the absence of
sinks or sources in the neighboring elements. The effectiveness of the PCA-step
will be backed up by the experiments in Section~\ref{sec:numericalExperiments}.

\subsection{Offline-Online Decomposition}
\label{sec:offlineOnlineDecomposition}
In this section, we will demonstrate how the computations for the LRBMS method
can be split into computationally demanding parts that will be executed during
the offline phase and computationally inexpensive parts to be performed during
the online phase. This is possible because of the affine splitting
(\ref{eq:parameterizedMob}).

Given a reduced basis $\redBas=\{\varphi_1,\ldots,\varphi_N\}$ of size $N\in\NN$, we compute all
parameter-independent parts of Equation (\ref{eq:coarseScalePressure}), that is:
all parts that depend on the location in space but not on $\Theta(s)$. We can
identify the terms $\underline{b}_q\in\RR^{N\times N},\, q=1,\ldots,\nParam$,
$\underline{c}\in\RR^{N\times N}$
\begin{align*}
(\underline{b}_q)_{ij} &= \sum_{\fineEl\in\fineGrid}\int_\fineEl \tMobPP{q} \perm \nabla \varphi_i \nabla \varphi_j \dx \\
	& \quad + \sum_{\fineIn\in\innerInt\cup\dirBounFP}\int_\fineIn
 \left( \mean{\tMobPP{q}\perm\nabla \varphi_i\cdot \n_\fineIn}\jump{\varphi_j} +  \mean{\tMobPP{q}\perm\nabla \varphi_j\cdot \n_\fineIn}\jump{\varphi_i} \right) \ds, \\
(\underline{c})_{ij} &= \sum_{\fineIn\in\innerInt\cup\dirBounFP} \frac{\sigma_\fineIn}{h_\fineIn} \int_\fineIn \jump{\varphi_i}\jump{\varphi_j} \ds,
\end{align*}
and the terms $\underline{d}_q,\,\underline{e}\in\RR^N$ for the right hand side:
\begin{align*}
(\underline{d}_q)_{i} &=
	\sum_{\fineEl\in\fineGrid}\int_{\fineEl} (\nwMobPP{q}\nwDensity+\wMobPP{q}\wDensity)\perm\,G\cdot\nabla \varphi_i\dx \\
	& - \sum_{\fineIn\in\innerInt\cup\dirBounFP}\int_\fineIn \mean{(\nwMobPP{q}\nwDensity+\wMobPP{q}\wDensity)\perm\,G\cdot\n_\fineIn}\jump{\varphi_i} \ds
	\quad- \sum_{\fineIn\in\dirBounFP}\int_{\fineIn} \left( \tMobPP{q}\perm\nabla \varphi_i\cdot \n_\fineIn \right) \dirValP \ds, \\
(\underline{e})_i &=
	\sum_{\fineEl\in\fineGrid}\int_{\fineEl} q_1 \varphi_i \dx
	+\sum_{\fineIn\in\dirBounFP}\int_{\fineIn}  \frac{\sigma_\fineIn}{h_\fineIn} \varphi_i  \dirValP \ds
	+ \sum_{\fineIn\in\neuBounFP}\int_\fineIn \neuValP \varphi_i \ds.
\end{align*}

During the online phase, for a given saturation $s\in\fineSpace$, the solution
$\underline{p}_H\in\RR^N$ of the discrete equivalent of the reduced equation
(\ref{eq:coarseScalePressure}) is then given by
\begin{align}
\label{eq:reducedSystem}
\left(\underline{c} + \sum_{q=1}^{\nParam} \theta_q(s) \underline{b}_q\right)\underline{p}_H = \underline{e}+\sum_{q=1}^{\nParam} \theta_q(s) \underline{d}_q.
\end{align}

While the computation of the quantities
$\underline{b}_q,\underline{c},\underline{d}_q$ and $\underline{e}$ has a
complexity polynomial in $1/h$, computing the sums in the
reduced system (\ref{eq:reducedSystem}) has a complexity polynomial in $N$. The
only critical part in (\ref{eq:reducedSystem}) is the computation of the
coefficients $\theta_i(s)$, $i=1,\ldots,\nParam$ that depends on the grid size
$h$ because of the least-squares approximation
(\ref{eq:leastSquaresFit}).
Nevertheless, as the complexity of the least-squares fit is only
$\mathcal{O}(1/h)$, we still expect largely accelerated computations compared
to a high-dimensional pressure solve, which has a complexity of $\mathcal{O}(1/h^2)$ or even $\mathcal{O}(1/h^3)$.
\begin{remark}
As mentioned before, different parametrizations like boundary value pa\-ram\-e\-trizations
would be possible for Problem (\ref{eq:mainPressure}-\ref{eq:mainBoundary}). As
long as these parametrizations are affine-linear, the offline-online
decomposition works analogously. For other kinds of parametrizations, see
e.g., \cite{Patera:un,Kaulmann2011rr,Haasdonk:2008io}.
\end{remark}


\section{LRBMS Two-Phase Flow Scheme} \label{sec:couplingSatAndPres} 
We now have all parts together to introduce our overall reduced approximation
scheme for two-phase flow in porous media using an IMPES-type coupling.
\begin{algorithm}[LRBMS Two-Phase Flow Scheme]
\label{alg:fullTwoPhaseScheme}
Compute the LRBMS pressures $\pC[1],\ldots,\pC[\nTimeSteps]$ and saturations
$\sC[1],\ldots,\sC[\nTimeSteps]$ as follows.
\begin{enumerate}
\item Choose a maximum basis size $N_{\max}\in\NN$, the number of mobility
profiles $\nParam\in\NN$ and the number of time steps $\nTimeSteps\in\NN$.
\item \label{itm:couplingBasisConstruction} Use Algorithm
\ref{alg:lrbmsBasisConstruction} to compute a reduced basis
$\redBas=\{\varphi_1,\ldots,\varphi_N\}\subset\fineSpace$ of size $N\in\NN$ for
mobility profiles $\wnwMobPP{1},\ldots,\wnwMobPP{\nParam}\in\fineSpace$.
\item Compute the quantities $\underline{b}_q, \underline{c}, \underline{d}_q$
and $\underline{e}$ from Section~\ref{sec:offlineOnlineDecomposition}.
\item Project the initial data $\sInitial$ to the fine grid:
	\begin{align*}
		\int_\Omega \sC[0]v\dx = \int_\Omega \sInitial v\dx,\quad\forall v\in\fineSpace.
	\end{align*}
\item\label{itm:LRBMSTwoPhaseSchemeTimeIteration} For $n=0,\ldots,\nTimeSteps-1$,
\begin{enumerate}
\item compute the reduced-dimensional pressure solution $\pD[n+1]\in\RR^N$ as
solution to Equation (\ref{eq:reducedSystem}) for $s=\sC[n]$. This means that we
need to compute the least-squares fit (\ref{eq:leastSquaresFit}) to the given
saturation $\sC[n]\in\fineSpace$ and then solve the reduced dimensional system
(\ref{eq:reducedSystem}).
\item \label{itm:couplingRecons} given the reduced-dimensional pressure solution
$\pD[n+1]$, reconstruct a function $p_h^r\in\fineSpace$ by setting
\begin{align*}
p_h^r = \sum_{i=1}^{N} \left(\pD[n+1]\right)_i \varphi_i.
\end{align*}
\item from the reconstructed pressure solution, compute the velocity
$\velF[n+1]\in\raviartThomas$ using (\ref{eq:velocityDG}) for $s=\sC[n]$ and
$p=p_h^r$.
\item compute the fine-scale saturation $\sC[n+1]\in\fineSpace$ using
(\ref{eq:saturationDG}).
\end{enumerate}
\end{enumerate}
\end{algorithm}
Notice the difference to the quantities computed in the high-dimensional
two-phase flow scheme: While in Algorithm~\ref{alg:highDimTwoPhaseFlow} we used
the linear mobilities (\ref{eq:brooksCorey}) in the DG-bilinear-form and in the
right hand side for the pressure, we now use the parametrized mobilities
(\ref{eq:parameterizedMob}). We use the notation $\sC[n]$ for saturations
computed with the LRBMS scheme---although the saturation itself is not computed
in a reduced space---to point out the dependency on the coarse-scale pressure.

The novelties in this scheme are the coupling between pressure and saturation
via a parametrized mobility and the replacement of the pressure equation by a
reduced-dimensional substitute. While both the least-squares approximation in
the basis construction step (\ref{itm:couplingBasisConstruction}) and the
reconstruction step (\ref{itm:couplingRecons}) still need to prove their
efficiency, we expect this scheme to allow largely accelerated computations with
acceptable additional error. The validity of this assumption will be
investigated in the numerical experiments. We therefore propose
Scheme~\ref{alg:fullTwoPhaseScheme} as an alternative to the reduction
approaches mentioned in Section~\ref{sec:introduction}.

\section{Numerical Experiments}
\label{sec:numericalExperiments}
In this section, we demonstrate the advantages of our approach introduced in
Sections~\ref{sec:lrbms} and \ref{sec:couplingSatAndPres} by means of a
2D-benchmark problem. All implementation was done using the Distributed and
Unified Numerics Environment (DUNE), see \cite{Bastian:2008fs, Bastian:2008dz,
Dedner:2010bq, Dedner:2012wq}. For the results shown in this section we will use
the true error in an energy norm as error measure:
\begin{align*}
\errorEst(\param, W)
	= \left[\pBilForm{\pC(\tMob)-\pF(\tMob)}
					 {\pC(\tMob)-\pF(\tMob)}
					 {\bar\lambda}
      \right]^{1/2},\quad\tMob = \sum \param_q\tMobPP{q}
\end{align*}
where $\bar\lambda=\sum_{q=1}^\nParam \bar{\param}_q \tMobPP{q}$ for a fixed
parameter $\bar{\param}$, $\pC(\tMob)$ denotes the LRBMS-pressure-solution in
the space $W$, see Definition~\ref{def:coarseScalePressure}, and $\pF(\tMob)$
denotes the high-dimensional pressure solution, see
Definition~\ref{def:pressureDG}. For the exposition of well-applicable
\textit{a-posteriori} estimators for our setting we refer to
\cite{Schindler:2014tz, Smetana:vx}.

\begin{figure}
	\centering
	\includegraphics[width=\textwidth]{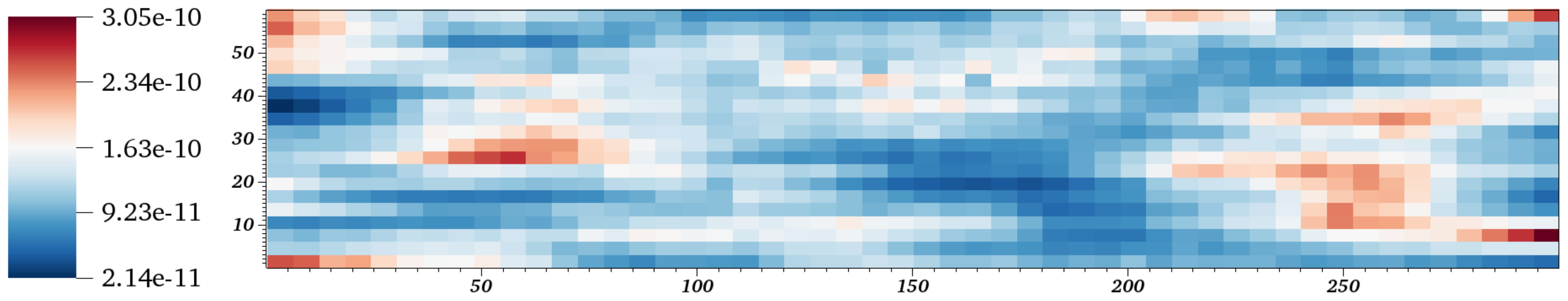}
	\caption{The permeability $\perm$ [m$^2$] used in the 2D benchmark problem}
	\label{fig:permEasy}
	\includegraphics[width=\textwidth]{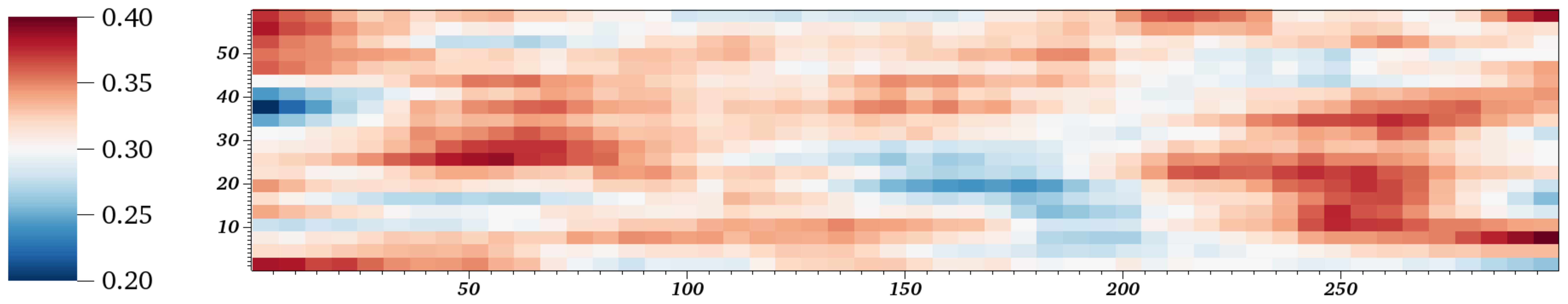}
	\caption{The porosity $\por$ used in the 2D benchmark problem}
	\label{fig:poroEasy}
\end{figure}

Our benchmark models the replacement of the non-wetting phase by the wetting phase in
$\Omega=[0,300]\times[0,60]$ for $T=3\cdot 10^5$.
The fine mesh $\fineGrid$ consists of $400 \cdot 160=64000$ rectangles and we
use $\nTimeSteps=6000$ time steps for the temporal discretization.
The densities are
$\wDensity=999.749$, $\nwDensity=890$, the viscosities are $\wVis=0.00130581$
and $\nwVis=0.008$. The permeability and porosity fields are shown in
Figures~\ref{fig:permEasy} and \ref{fig:poroEasy}, respectively.

\begin{figure}
\begin{minipage}{0.49\textwidth}
\centering
\includegraphics[width=\textwidth]{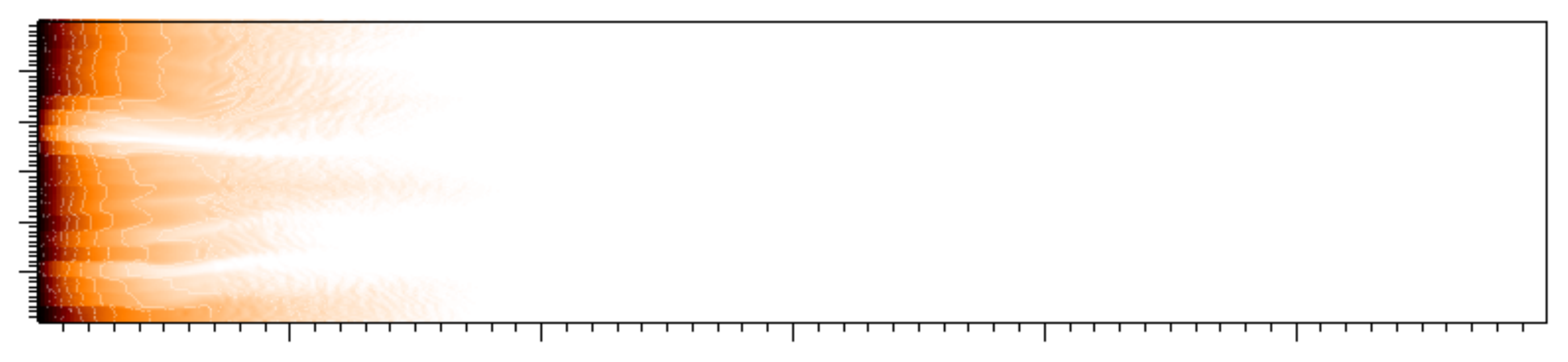} \\
\includegraphics[width=\textwidth]{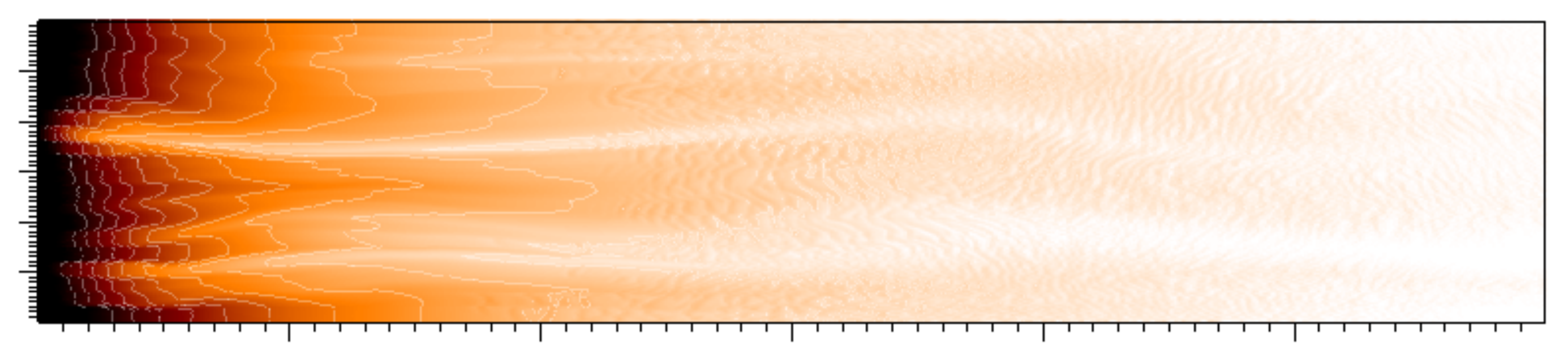}
\end{minipage}
\begin{minipage}{0.49\textwidth}
\centering
\includegraphics[width=\textwidth]{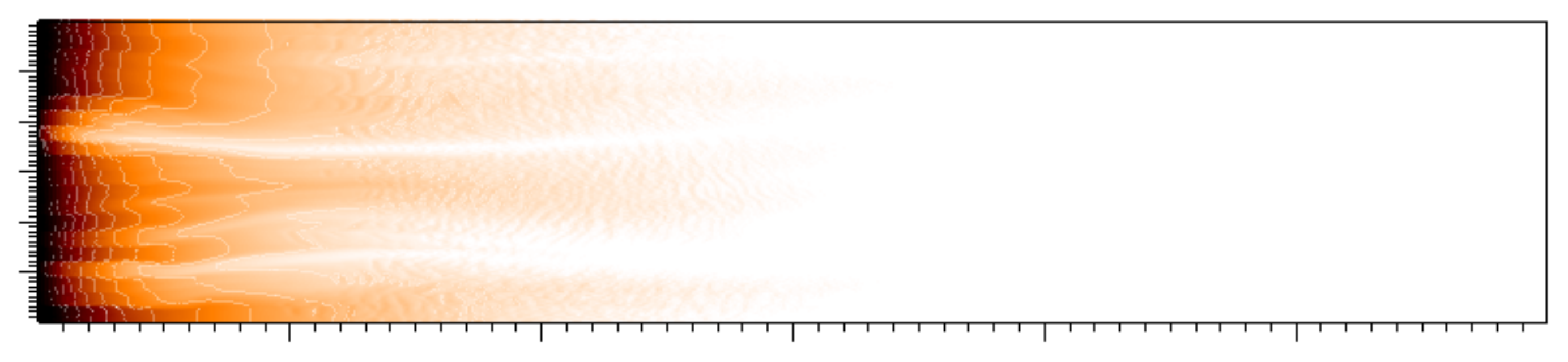} \\
\includegraphics[width=\textwidth]{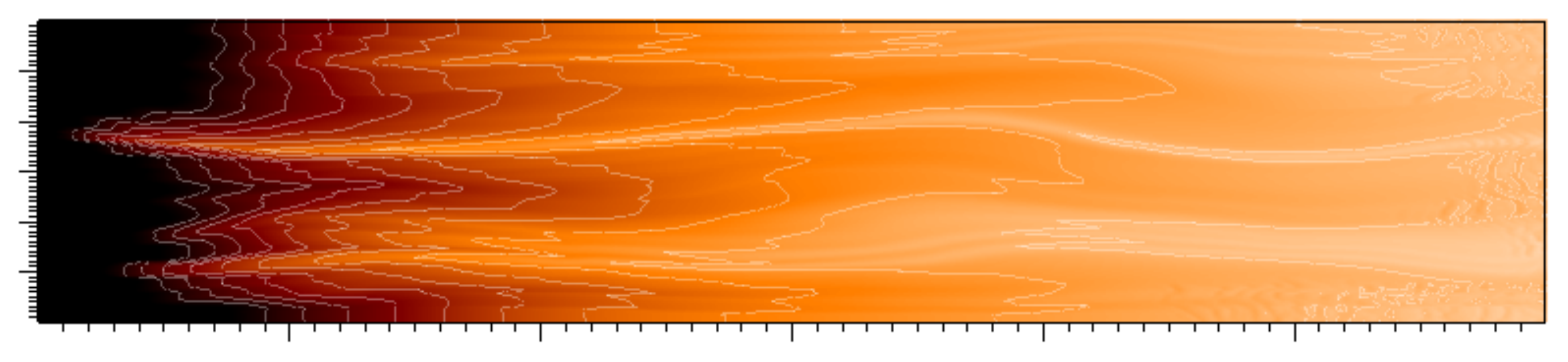}
\end{minipage}
\begin{center}
\includegraphics[width=0.5\textwidth]{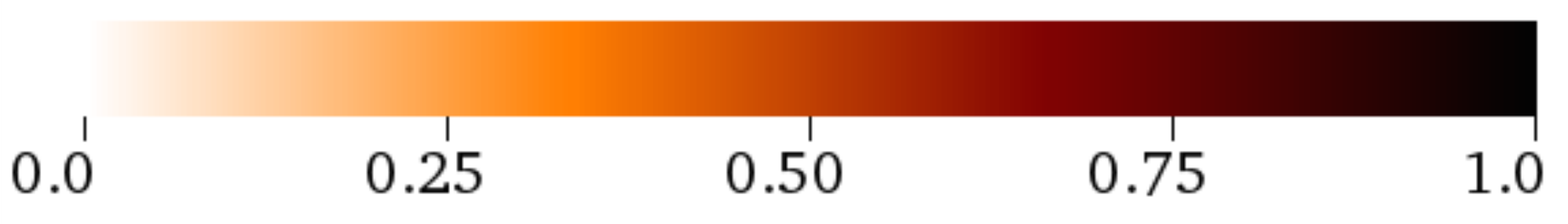}
\end{center}
\caption{Wetting phase saturation $\sF$ computed using the LRBMS method after 3.5, 10,
15 and 48 hours (from left to right and top to bottom) including contour lines}
\label{fig:saturationSnapshots}
\end{figure}

The domain is initially fully saturated with the non-wetting phase
($\sInitial\equiv0$), which is then displaced by the wetting phase entering from the
left boundary, modelled via the boundary conditions
\begin{align*}
\begin{aligned}
	s &= 1.0 & \text{in } \dirBounFS\times [0,T], \\
	p &= 10 & \text{in } \dirBounFP\times [0,T],
\end{aligned}
\qquad
\begin{aligned}
 - \tMobBC(s) \perm \nabla p\cdot\n &= 3\cdot10^{-4} & \text{in } \neuBounFP^1\times [0,T], \\
 - \tMobBC(s) \perm \nabla p\cdot\n &= 0 & \text{in } \neuBounFP^2\times [0,T],
\end{aligned}
\end{align*}
on $\dirBounFS=\dirBounFP=\{0\}\times[0,60]$ and
$\neuBounFP=\neuBounFP^1\cup\neuBounFP^2$, where
$\neuBounFP^1=\{300\}\times[0,60]$,
$\neuBounFP^2=[0,300]\times\{0\}\cup[0,300]\times\{60\}$. In this benchmark no
sources are used ($q_1\equiv q_2\equiv 0.0$) and we neglect gravity so that
$\gravity=(0.0,0.0)^\transpose$. The relative permeabilities in Equation
(\ref{eq:brooksCorey}) are given here via the linear relations
$k_\text{rw}(s)=s$ and $k_\text{ro}(s)=1-s$ for simplicity. Different relations
like the Brooks-Corey or van Genuchten law would be possible.
Figure~\ref{fig:saturationSnapshots} shows the saturation $\sF$ computed with
the full scheme (Algorithm~\ref{alg:highDimTwoPhaseFlow}) after approximately
3.5, 10, 15, and 48 hours.

\begin{figure}
	\centering
	\includegraphics[width=\textwidth]{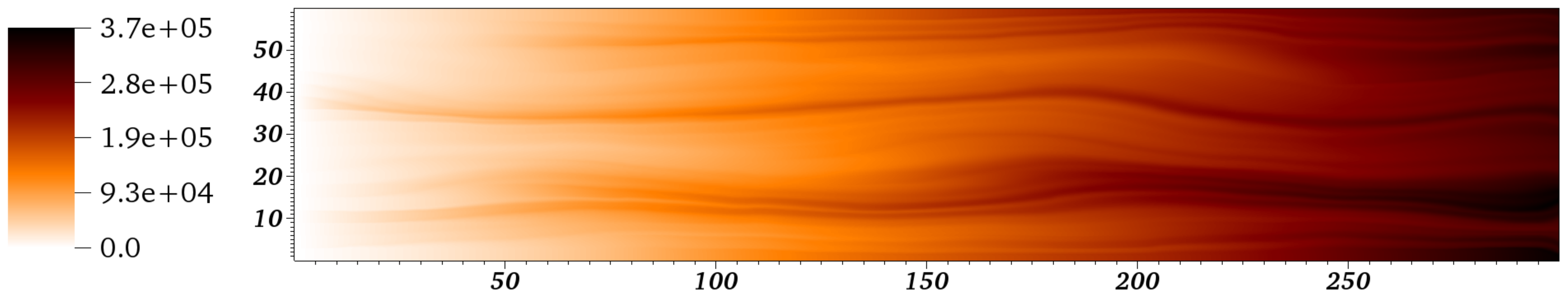}
	\caption{The time-of-flight $\tofF{s}$ for the 2D benchmark problem for $\s=0$}
	\label{fig:tofEasy}
\end{figure}

\begin{figure}
	\begin{minipage}{0.49\textwidth}
		\centering
		\includegraphics[width=\textwidth]{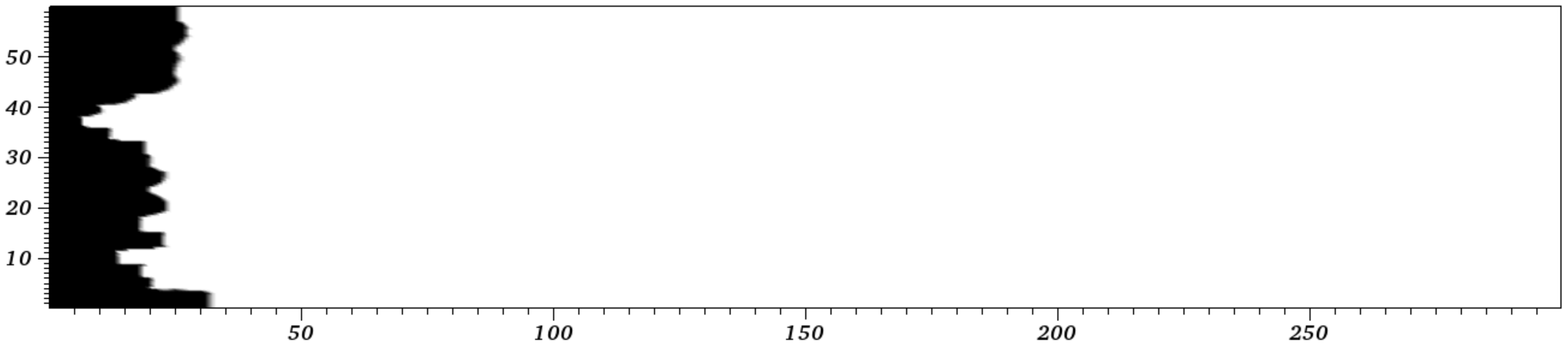}
		\includegraphics[width=\textwidth]{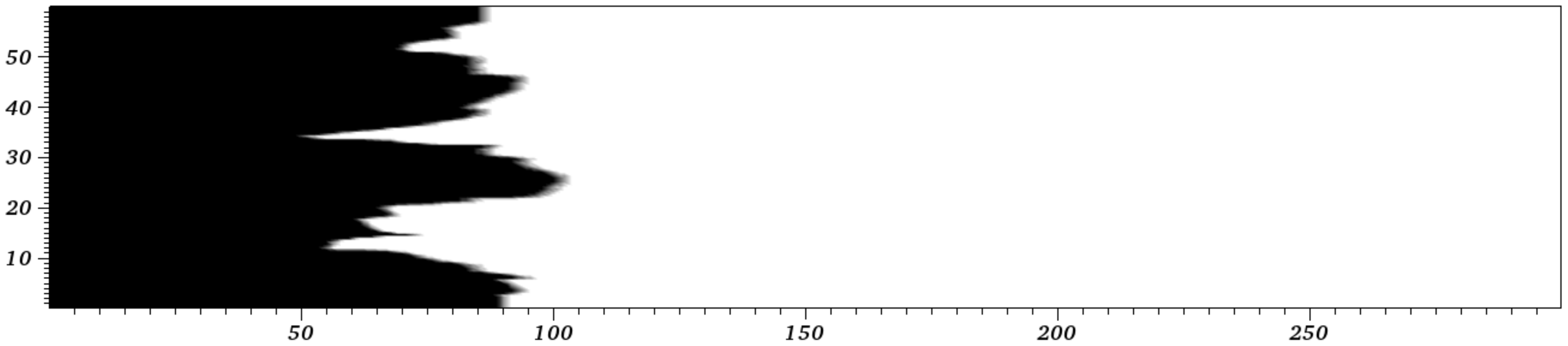}
		\includegraphics[width=\textwidth]{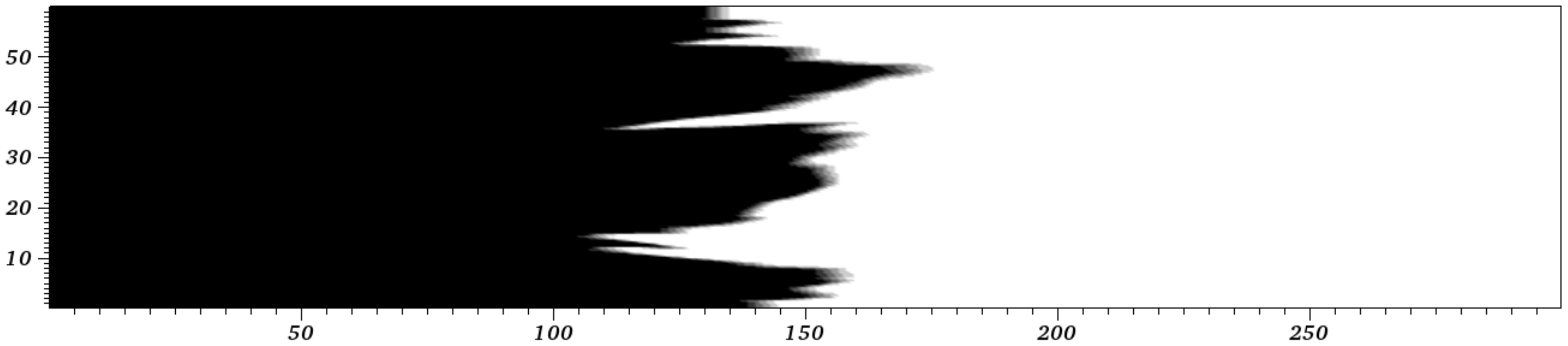}
	\end{minipage}
	\begin{minipage}{0.49\textwidth}
		\centering
		\includegraphics[width=\textwidth]{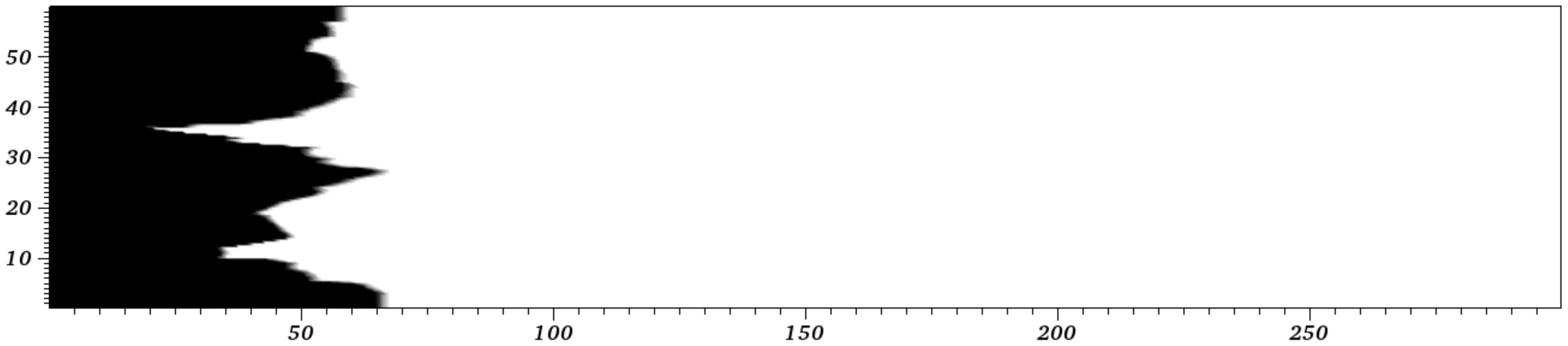}
		\includegraphics[width=\textwidth]{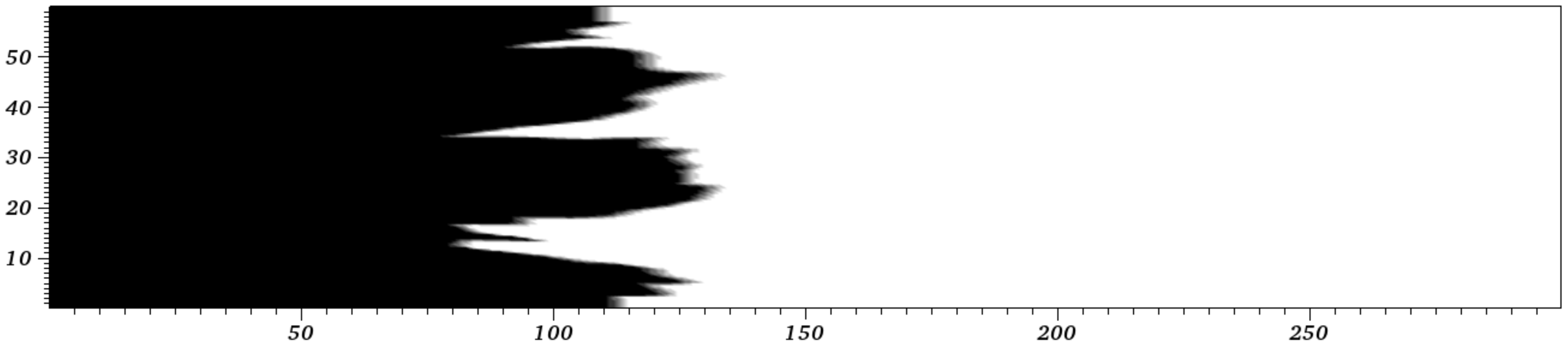}
		\includegraphics[width=\textwidth]{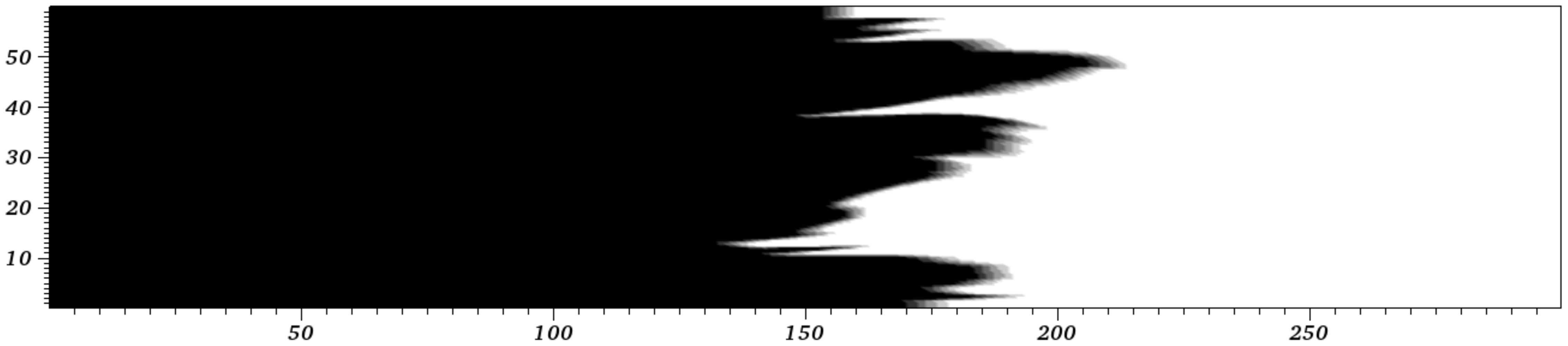}
	\end{minipage}
	\begin{center}
	\includegraphics[width=0.4\textwidth]{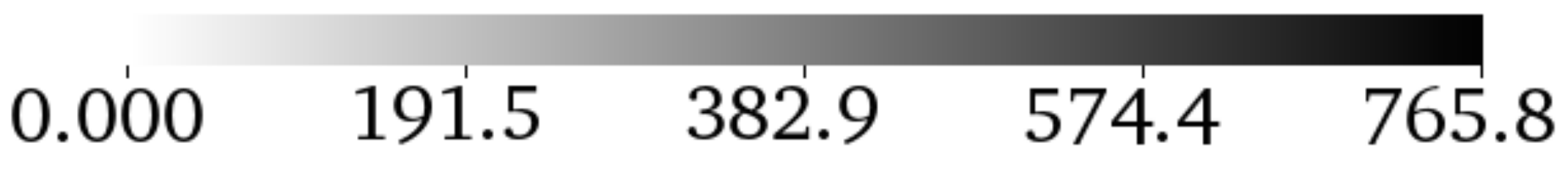}
	\end{center}
	\caption{Wetting mobility profiles computed with
	Algorithm~\ref{alg:lrbmsBasisConstruction} for the time-of-flight depicted in
	Figure~\ref{fig:tofEasy}. Not shown are the profiles for $q=1$ and $q=M$ which
	have constant values $\wMobPP{0}\equiv0$ and $\wMobPP{M}\equiv765.808$,
	respectively.}
	\label{fig:mobilitiesEasy}
\end{figure}

From the time-of-flight, which is depicted in Figure~\ref{fig:tofEasy}, we
compute eight wetting and non-wetting mobility profiles using
Algorithm~\ref{alg:lrbmsBasisConstruction}. The resulting wetting mobility
profiles are depicted in Figure~\ref{fig:mobilitiesEasy}.
Using the profiles we compute the reduced basis $\redBas$ using coarse meshes
with sizes $|\coarseGrid|=1,4,8,16,32$.
We use the tolerance $\epsilon_\text{tol}=10^{-4}$, a
training set $\trainSet$ consisting of $300$ randomly distributed parameters
$\param$ in $[0.0001, 1]^8$ with $\sum(\param)_i=1$, and the maximum size
$N_{\max}=\infty$ for Algorithm~\ref{alg:lrbmsBasisConstruction}. The resulting
basis sizes can be seen in Table~\ref{tab:basisSizesEasyFlow}. We observe that
with increasing size of the coarse mesh (first column), the number of snapshots
computed in Step~\ref{itm:lrbmsGreedySnapshot} of
Algorithm~\ref{alg:lrbmsBasisConstruction} (fourth column) decreases
significantly from 134 to 69. At the same time, the overall basis size (the sum
of all local basis sizes, second column) increases. Notice that the basis size
is not necessarily equal the product between the number of snapshots and the
coarse grid size. This is because we orthonormalize each new snapshot with
respect to the existing basis in each extension during the basis generation,
reject local extensions with norms below a certain threshold to avoid linear
dependencies in the resulting basis and hence reduce the number of local basis
functions.

\begin{table}
\centering
\includegraphics[width=0.5\textwidth, trim=0cm -0.5cm 0cm 0cm]{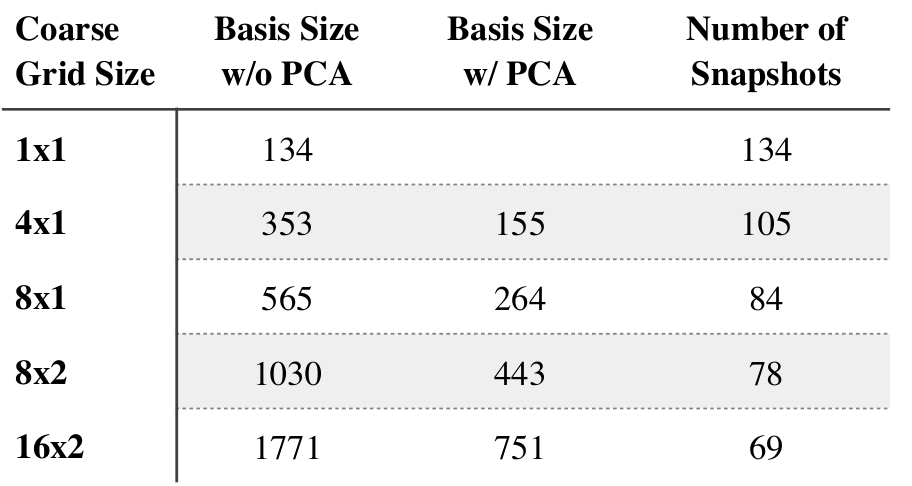}
\caption{Basis sizes $|\redBas|$ resulting from
Algorithm~\ref{alg:lrbmsBasisConstruction} for the 2D benchmark problem and a
fine triangulation with 64000 elements: Size of the coarse mesh, sum of all
local basis sizes before and after application of the PCA, number of snapshots
computed during the basis generation}
\label{tab:basisSizesEasyFlow}
\end{table}

In Table~\ref{tab:basisSizesEasyFlow} we also see the impact of the local data
compression using PCA (Step~\ref{itm:lrbmsGreedyPCA} in
Algorithm~\ref{alg:lrbmsBasisConstruction}): With increasing coarse mesh size,
the PCA is able to reduce the local basis sizes significantly on the different
coarse elements.

\begin{figure}
\centering
\includegraphics[width=0.75\textwidth]{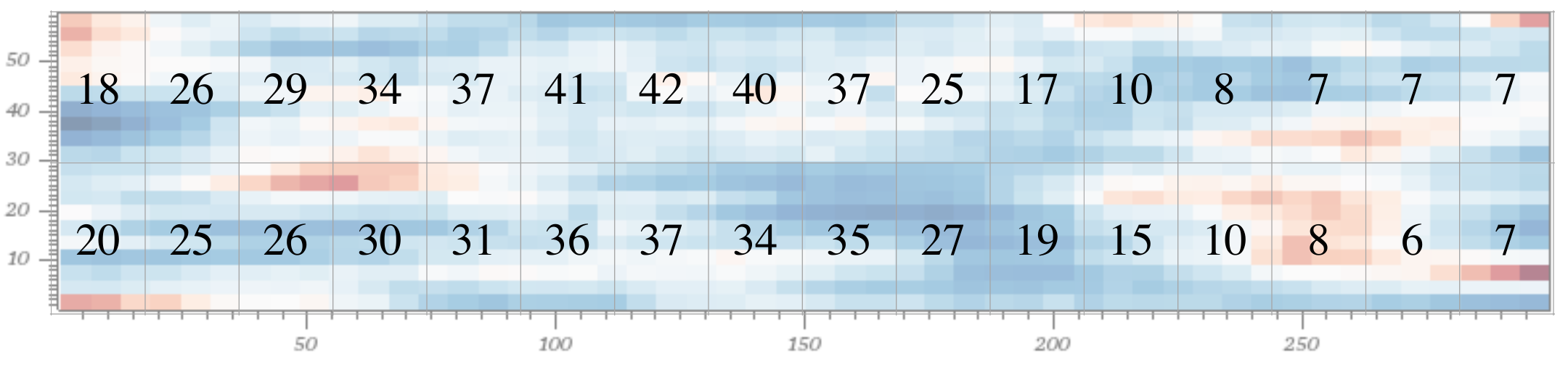}
\caption{Foreground: local basis sizes after the basis generation performed for
Table~\ref{tab:basisSizesEasyFlow}, last line, including the PCA. Background:
permeability field}
\label{fig:impactOfPCA}
\end{figure}

This effect can be seen again in Figure~\ref{fig:impactOfPCA} where we plot the
local basis sizes after the PCA step for the test run presented in
Table~\ref{tab:basisSizesEasyFlow}, last line. We observe that the local basis
size varies strongly from left to right, due to the fact that the peak
saturation does not reach the right half of the computational domain before the
end time $T$. Also, the basis size differs from bottom to top, especially in
regions where the permeability, which is plotted in the
background of Figure~\ref{fig:impactOfPCA}, shows strong variation from bottom
to top.

In Tables~\ref{tab:dataSizesEasySat} and \ref{tab:dataSizesEasyPress} we see the
resulting discrepancies between the saturation $\sF[n]$  computed by the SWIP-DG
method on the fine mesh (Algorithm~\ref{alg:highDimTwoPhaseFlow}) and the
saturation $\sC[n]$ computed with the reduced scheme
\ref{alg:fullTwoPhaseScheme}, as well as the discrepancy between the SWIP-DG
pressure and the reduced pressure during the two-phase flow simulation. We
present both the relative $\lTwo$ and $\hOne$ discrepancies for the saturation
\begin{align*}
\lerror^n=\frac{\| \sF[n]-\sC[n]\|_{\lTwo(\Omega)}}{\|\sF[n]\|_{\lTwo(\Omega)}},\quad
\herror^n=\frac{\| \sF[n]-\sC[n]\|_{\hOne(\Omega)}}{\|\sF[n]\|_{\hOne(\Omega)}},
\end{align*}
and the respective quantities for the pressure. For different coarse grids, the
second column in each table displays the number of snapshots computed during the
extension step \ref{itm:lrbmsGreedySnapshot} in
Algorithm~\ref{alg:lrbmsBasisConstruction} as a measure for the amount of work
needed in that step. Then, in the two sections for the $\lTwo$ and $\hOne$
relative discrepancies, we present mean discrepancies over all time steps and
the discrepancy at the last time step $\tEnd$ for basis generations without and
with usage of the PCA.

Over the different coarse mesh configurations, we consistently observe a mean
$\lTwo$-discrepancy for the saturation of approximately $5.6$\% (standard
deviation: $1.6$\%), both with and without the PCA, which is reduced to $4.3$\%
at end time. The respective $\hOne$-discrepancies are slightly bigger with a
mean of $7.9$\% (standard deviation: $2.9$\%) and $5.5$\% at end time. The mean
$\lTwo$-discrepancy for the pressure is approximately $1.4$\% (standard
deviation: $0.3$\%) for all coarse grid configurations, and reduces to $1.3$\%
at end time. The respective $\hOne$-discrepancies are in the same ranges.
In conclusion, we can say that the data compression step using the PCA increases
neither the error for the pressure nor the error for the saturation noticeably.

\begin{figure}
\centering
\includegraphics[width=\textwidth]{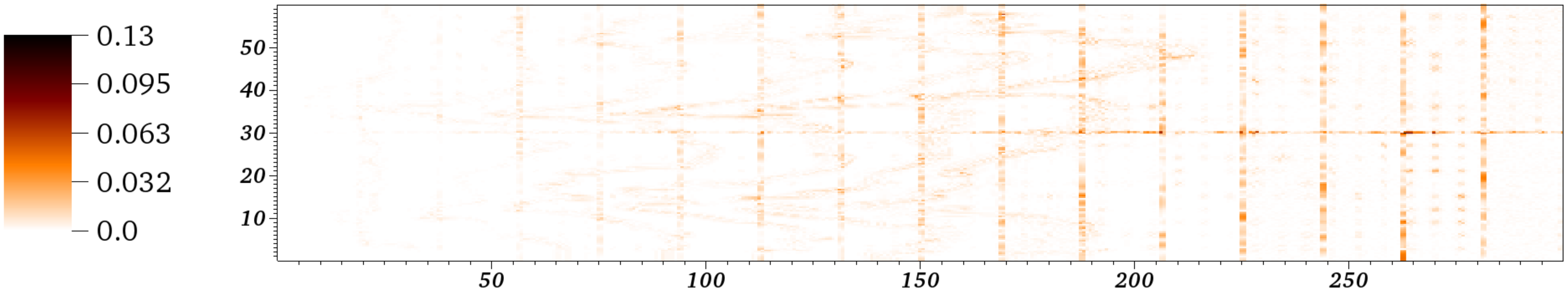}
\caption{Relative mass loss $\zeta^n$ for the final time step $t^n=3\cdot 10^5$}
\label{fig:easyHomConservativity}
\end{figure}

While we can ensure mass conservation on the coarse grid by adding
per-coarse-cell unit basis functions to the reduced basis, the method is not
mass-conservative on the fine grid. In Figure~\ref{fig:easyHomConservativity} we
plot the relative mass loss
\begin{align*}
\zeta^n:\fineGrid\to\RR^+,\quad \zeta^n(\fineEl)=
	\frac{1}{\max_{x\in\partial\fineEl}\|\velF[n](x)\|_{l^2}}
	\left|
		\int_{\partial\fineEl} \velF[n]\cdot\n \ds
	\right|
\end{align*}
as a measure of the lack of mass conservation at end time $t^n=3\cdot 10^5$ for
the velocity computed with the LRBMS scheme on a $16\times 2$ coarse grid. We
see that the velocity lacks mass conservation mainly on coarse cell
intersections and in regions with large gradients in the mobility profiles
$\tMobPP{i}$. Overall, the relative lack of mass conservation is well below
$2\cdot 10^{-5}$.
Although we therefore found this to be only a minor problem in the tests---the
saturation usually did not grow above a value of $1.05$ and hardly ever fell
below zero---future work will include research on how to ensure
mass-conservation on the fine grid.

As the discrepancies are increased neither by introducing more coarse cells, nor
by application of the PCA, we find the Localized Reduced Basis Multiscale method
to be well applicable in this context: Introduction of more coarse cells reduces
the number of costly snapshots to be computed during the basis construction,
while application of the PCA yields smaller, more compact bases, leading to
faster computations during the two-phase flow simulation.

\begin{table}
	\centering
	\includegraphics[width=\textwidth, trim=0cm -0.5cm 0cm 0cm]{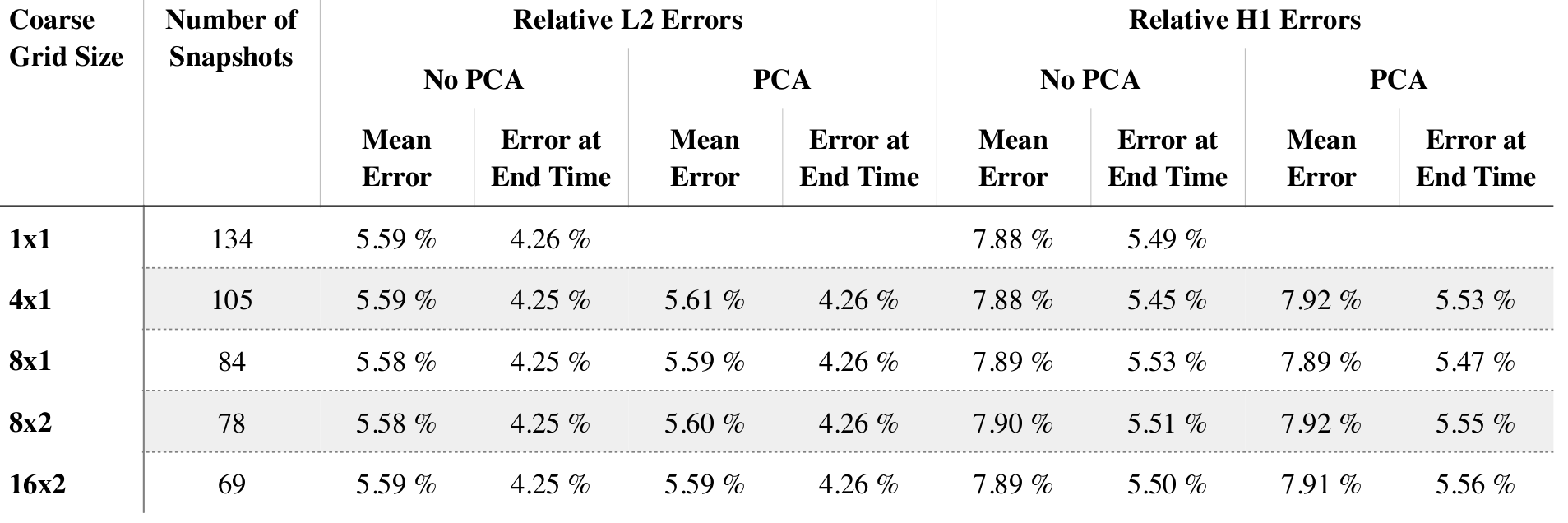}
	\caption{Relative $\lTwo$ and $\hOne$ discrepancies between the saturation
	computed with Algorithm~\ref{alg:fullTwoPhaseScheme} and the saturation
	computed with the fine-scale algorithm \ref{alg:highDimTwoPhaseFlow} for
	different coarse mesh configurations: Number of snapshots needed during the
	basis construction, mean relative $\lTwo$ discrepancy and relative $\lTwo$
	discrepancy at end time $T$ without using the PCA, respective quantities after
	usage of the PCA, respective quantities for the $\hOne$ norm}
	\label{tab:dataSizesEasySat}
\end{table}

\begin{table}
	\centering
	\includegraphics[width=\textwidth, trim=0cm -0.5cm 0cm 0cm]{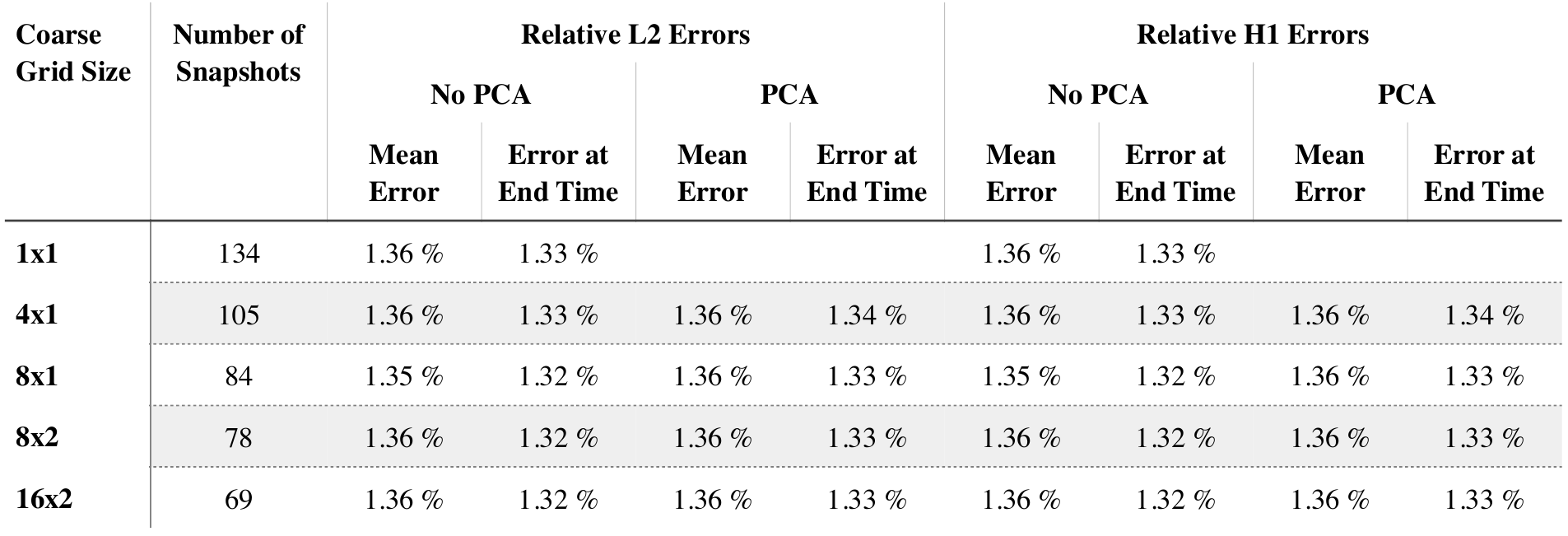}
	\caption{Relative $\lTwo$ and $\hOne$ discrepancies between the pressure
	computed with Algorithm~\ref{alg:fullTwoPhaseScheme} and the pressure computed
	with the fine-scale algorithm \ref{alg:highDimTwoPhaseFlow} for different
	coarse mesh configurations: Number of snapshots needed during the basis
	construction, mean relative $\lTwo$ discrepancy and relative $\lTwo$
	discrepancy at end time $T$ without using the PCA, respective quantities after
	usage of the PCA, respective quantities for the $\hOne$ norm}
	\label{tab:dataSizesEasyPress}
\end{table}

Test runs with a different tolerance $\epsilon_\text{tol}=10^{-5}$ for the
Greedy basis construction (Algorithm~\ref{alg:lrbmsBasisConstruction}) exhibit
the same discrepancies in both the $\lTwo$ and $\hOne$ norms. Further,
computations with higher numbers of mobility profiles ($\nParam=10, 12, 20$)
give roughly the same discrepancies. This gives rise to the assumption that the
error is dominated by different phenomena, which we consider to be twofold:
First, the assumption that the time-of-flight is invariant for the whole simulation
is not valid in some regions. Therefore, the position of the saturation front
and its impact on the mobility cannot be represented correctly. Second, the
profile of the mobility along streamlines is not approximated well enough due to
the ad-hoc profile generation in Algorithm~\ref{alg:lrbmsBasisConstruction}.

\begin{figure}
\centering
\includegraphics[width=\textwidth]{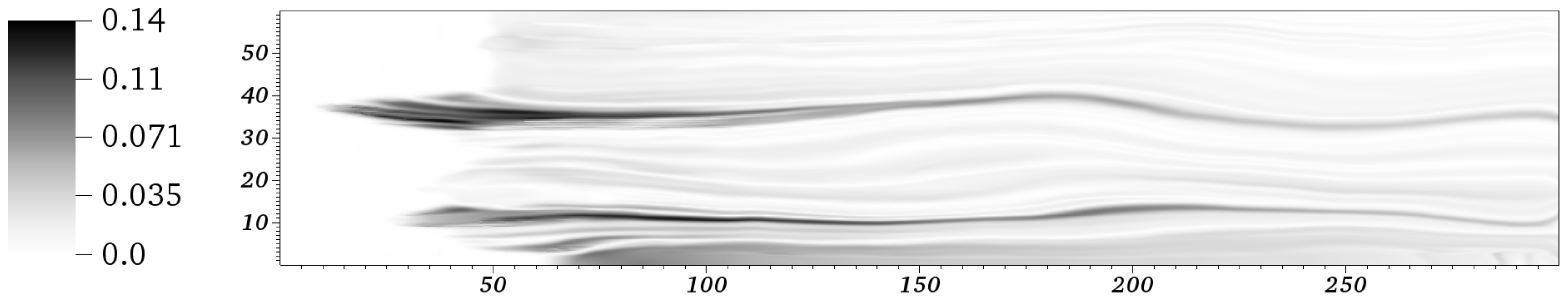}
\caption{Absolute value of the difference in the saturation $\sF[\nTimeSteps]$ computed with the full high-dimensional scheme and the saturation $\sC[\nTimeSteps]$ computed with the LRBMS scheme at end time $T$}
\label{fig:saturationDifference}
\end{figure}

\begin{figure}
\centering
\includegraphics[width=\textwidth]{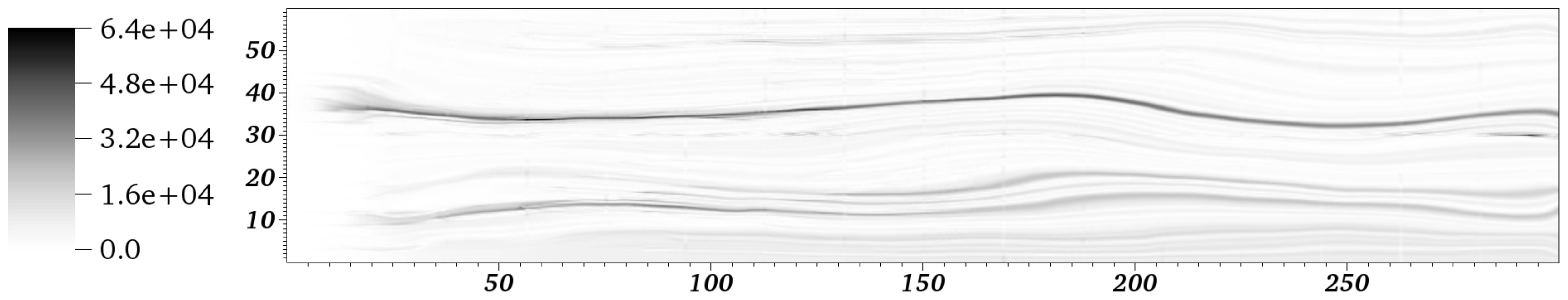}
\caption{Absolute difference between the time-of-flight for the initial
saturation $\sInitial$ and the time-of-flight for the saturation computed with
the LRBMS scheme for a coarse grid with $16\times2$ cells at time $t=10000s$}
\label{fig:tofDifference}
\end{figure}

To support the first statement, we present in
Figure~\ref{fig:saturationDifference} the absolute value of the difference
between the SWIP-DG approximation $\sF[\nTimeSteps]$ and our approximation of
the saturation $\sC[\nTimeSteps]$ at end time $\tEnd$. We see that huge
discrepancies arise in three distinct positions: around the point
$\mathbf{x}=(30, 35)$, the point $\mathbf{x}=(40, 10)$, and the point
$\mathbf{x}=(80,1)$. These positions are located directly downwind along the
streamlines from points where the time-of-flight changes significantly during
the test run, see Figure~\ref{fig:tofDifference}. The error produced in those
regions is then transported through the domain along the streamlines, hence the
error distribution to be seen in Figure~\ref{fig:saturationDifference} is
established.

The second statement can be justified by using $\nParam$ high-dimensional
saturation approximations $\sF[n_1],\ldots,\sF[n_\nParam]$ at points
$t^1,\ldots,t^\nParam$ in time to form the profiles in
Step~\ref{itm:lrbmsBasisConsInitEnd} of
Algorithm~\ref{alg:lrbmsBasisConstruction}:
\begin{align*}
\wMobPP{q}=\wMobBC(\sF[q]),\quad \nwMobPP{q}=\nwMobBC(\sF[q]).
\end{align*}
In doing so, we ensure that the shape of the mobility profiles along streamlines
is correct and the error of the LRBMS two-phase flow scheme should decrease
drastically. Indeed, for $\nParam=8$, this procedure decreases the mean relative
$\lTwo$-discrepancy to 2.3\% with a standard deviation of 2.0\%. Even more: The
$\lTwo$-discrepancy at $t=T$ is reduced by a factor of ten to 0.45\%. Obviously,
using multiple high-dimensional saturation profiles is not possible in general as the
model order reduction would become superfluous, but two other remedies could be
implemented: One could be to recompute the mobility profiles and the reduced
basis after a certain number of time steps $t^n$ using the time-of-flight for
$\sF[n]$. Another approach could be to solve the full high-dimensional two-phase
flow problem along some streamlines in one dimension in a preparatory step and
use a combination of the resulting mobility profiles as a replacement for those
in Step~\ref{itm:lrbmsBasisConsInitEnd} of
Algorithm~\ref{alg:lrbmsBasisConstruction}.


In Table~\ref{tab:runTimesOffline} we present runtimes for the basis
construction part of the LRBMS two-phase flow scheme. Again, we give the coarse
grid configuration (column one) and the number of snapshots needed (column two).
The third column then shows the time needed for the basis construction
(Algorithm~\ref{alg:lrbmsBasisConstruction}). The timings presented here include
the time for the computation of the mobility profiles (about 13 seconds), the
so-called ``training''-step (computation of all reduced solutions, evaluation of
the error measure, selection of the parameter for basis extension), and the
extension-step including the computation of a high-dimensional snapshot and
application of the Gram-Schmidt procedure. It does not contain the application
of the PCA which is consistently below one second for all coarse grid
configurations and hence can be considered negilgible. The time needed to
compute the reduced basis $\redBas$ is 50 minutes on one coarse cell (which
corresponds to a standard Reduced Basis method), goes up to 1 hour 7 minutes for
a coarse grid of size $4\times 1$ and is then reduced to 45 minutes for all
other coarse grids. The increase in total runtime from line one to line two can
explained by the relatively high number of iterations that is still needed to
reach the error tolerance with an increased per-step cost due to the larger
reduced bases. This effect would vanish if the snapshot computation (column
four) was more expensive, for larger fine grids, for example. As mentioned
earlier, we use the true error as error measure. The time for error estimation
in column five includes the time needed to  reconstruct a high-dimensional
snapshot from each reduced solution and compute the difference in the
$\lTwo$-norm. The time for computation of the high-dimensional snapshots itself
is not included.

In conclusion we can say that for the test case at hand, introduction of more
coarse cells yields slight reductions in terms of runtime. As the computation of
the snapshots shows a speedup by a factor of two for the computation on 32
coarse cells compared to the computation on one coarse cell, we can expect the
runtime-gain to increase drastically as the size of the fine mesh---and hence
the time for the computation of one high-dimensional snapshot---is increased.

\begin{table}
\centering
\includegraphics[width=0.6\textwidth, trim=0cm -0.5cm 0cm 0cm]{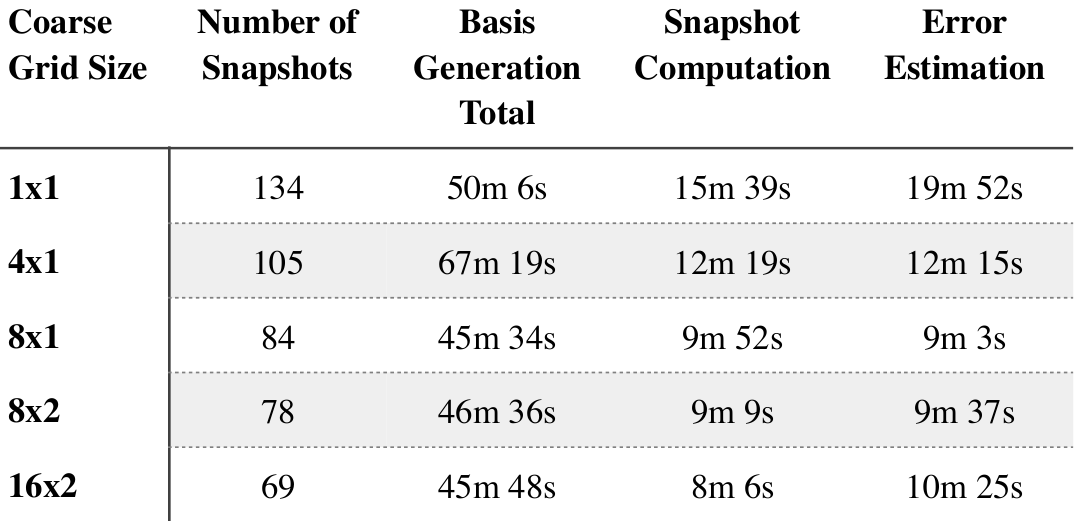}
\caption{Runtimes for the basis construction
(Algorithm~\ref{alg:lrbmsBasisConstruction}): Coarse grid size, number of
snapshots computed, total runtime, total time for snapshot computation, total
time for error estimation}
\label{tab:runTimesOffline}
\end{table}

Finally, in Table~\ref{tab:runTimesOnline} we present runtimes for the two-phase
flow simulation (Step~\ref{itm:LRBMSTwoPhaseSchemeTimeIteration} in
Algorithm~\ref{alg:fullTwoPhaseScheme}) for uncompressed bases (that is: bases
that were computed without usage of the PCA) and compressed bases on different
coarse grids and, for comparison, the same runtimes for a full high-dimensional
computation. We see that for the reduced simulations, more than 50\% of the
overall runtime of about 2 hours 40 minutes is spend in the application of the
ODE solver and slope limiter for the saturation equation. About one hour is
spend for the elliptic equation with about 45 minutes for the reconstruction of
the flux (see Definition~\ref{def:velocityDG}). For uncompressed bases,
computing all pressure solutions takes five to 15 minutes, depending on the
coarse grid and respective basis size, for compressed bases those times drop to
two to five minutes.

The time to solution for one reduced pressure computation therefore ranges from 20
milliseconds on the $1\times1$, $4\times1$ and $8\times1$ coarse grids to 50
milliseconds on the $8\times2$ and $16\times2$ coarse grids using the PCA.
Comparing these runtimes to the time needed for a high-dimensional simulation we
see the advantage of our method: The high-dimensional test run takes
approximately 16 hours with nearly 90\% of the time spent in the treatment of
the elliptic equation: more than two hours are spent assembling the pressure
system, solving it takes approximatly 11 hours in total (approximately seven
seconds per solve). The speed-up for the solution of the pressure equation is
approximately a factor 140. Because nothing was done to speed up the transport
solve, the overall speed-up for the two-phase flow simulation (high-dimensional
vs. basis construction and reduced simulation) is five. To also speed up the
transport solve significantly, one could replace the explicit temporal
discretization by a backward Euler scheme and utilize the fact that the
resulting nonlinear system share the same unidirectional flow properties as the
time-of-flight equation and hence can be computed in a per-element fashion with
local control over the nonlinear iterations; see \cite{NatvigLie:2008:jcp} for
details. Note that in this case, the slope limiter presented in Section
\ref{sec:highDimDisc} needs to be replaced by a different kind of stabilization
as it is not compliant with implicit time stepping schemes.

\begin{table}
\centering
\includegraphics[width=\textwidth, trim=0cm -0.5cm 0cm 0cm]{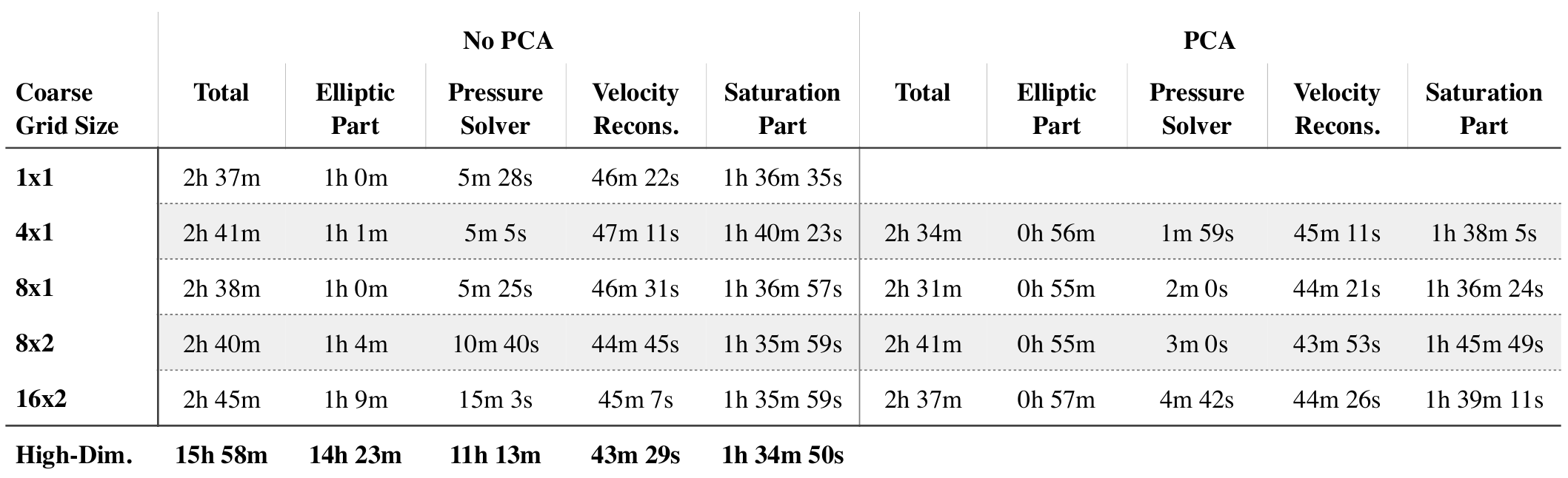}
\caption{Runtimes for Step~\ref{itm:LRBMSTwoPhaseSchemeTimeIteration} of
Algorithm~\ref{alg:fullTwoPhaseScheme} for uncompressed and compressed bases:
Coarse grid size, total runtime, time spent for treatment of the elliptic equation, thereof time spent in pressure
solve and velocity reconstruction and time spent in saturation part. Last line:
same numbers for a full high-dimensional simulation}
\label{tab:runTimesOnline}
\end{table}


\section{Summary and Outlook}
\label{sec:summaryAndOutlook}

We introduced a localized version of the Reduced Basis method. The idea of our
approach is to make use of two computational meshes: A fine mesh is used to
compute detailed solutions of a parametrized elliptic equation for different
parameters. Those solutions are then localized to the cells of a second, coarser
mesh. An optional data compression is applied and the result is used as basis
for a reduced-dimensional surrogate of the high-dimensional scheme. We applied
this technique to the pressure equation in a two-phase flow setting replacing
the orignal mobility by a parametrized approximation. We were able to
demonstrate significant reduction of the computational effort at acceptable
error compared to a full detailed two-phase flow simulation.

The quality of the reduced approximations for saturation and pressure is mainly
decided by the quality of the approximation of the original mobility by our
parametrized surrogate, therefore future work will include improvements of the
mobility approximation.

Mass conservation can only be guaranteed on the coarse mesh of our scheme, in
future work we will address the problem of mass conservation on the fine mesh
for the reduced simulations.
Applying the LRBMS to the saturation equation also would be an extension promising additional speedup.


\section*{Acknowledgement}
The author S.~Kaulmann would like to thank the German Research Foundation
(DFG) for financial support of the project within the Cluster of Excellence in
Simulation Technology (EXC 310/1) at the University of Stuttgart and within the
German Priority Programme 1648, ``SPPEXA - Software for Exascale Computing''.
Furthermore, financal support by the Baden-Württemberg Stiftung gGmbH is
gratefully acknowledged. Parts of the work was carried out while the first
author visited SINTEF in Oslo, Norway.



\end{document}